\documentclass[12pt]{article}
\usepackage[latin1]{inputenc}
\usepackage[T1]{fontenc}
\usepackage{fixmath}
\usepackage{amsmath}
\usepackage{amssymb}
\usepackage{natbib}
\usepackage{graphicx}
\usepackage{subfigure}
\usepackage{xcolor}
\usepackage{multirow}
\usepackage{mathabx}

\newtheorem{rem}{Remark}[section]
\newtheorem{lemma}{Lemma}[section]
\newtheorem{theo}{Theorem}[section]
\newtheorem{cor}{Corollary}[section]

\def\nPto{\xrightarrow[n\to\infty]{P}}
\def\beq{\begin{eqnarray*}}
\def\eeq{\end{eqnarray*}}
\def\e{\varepsilon}
\def\Var{\text{Var}}
\def\boxi{\hfill$\Box$}

\def\bx{{\mathbf x}}
\def\bX{{\mathbf X}}
\def\bo{{\boldsymbol \omega}}
\def\bbo{\bar{\boldsymbol \omega}}
\def\bbX{\bar{{\mathbf X}}}
\def\bbx{\bar{{\mathbf x}}}
\def\f{\text{\rm fun}}
\def\c{\text{\rm cat}}
\def\cS{\mathcal S}

\title{Uniform convergence rates and automatic variable selection in nonparametric regression with functional and categorical covariates}

\author{Leonie Selk\footnote{Helmut-Schmidt-University, Department of Mathematics and Statistics, Hamburg, Germany, corresponding author: leonie.selk@uni-hamburg.de}}

\begin{document}

\maketitle

\begin{abstract}
In \citet{SelkGertheiss2022} a nonparametric prediction method for models with multiple functional and categorical covariates is introduced. The dependent variable can be categorical (binary or multi-class) or continuous, so that both classification and regression problems are considered. In the paper at hand the asymptotic properties of this method are studied. A uniform convergence rate for the regression / classification estimator is given. It is further shown that a data-driven least squares cross-validation method can asymptotically remove irrelevant noise variables automatically.
\end{abstract}


 
\section{Introduction}
This paper deals with nonparametric prediction and estimation with multiple predictors that admit a mix of functional and categorical ones. The method considered was introduced in \citet{SelkGertheiss2022}. 
The idea is based on the well-known Nadaraya-Watson estimator 
\[\hat{m}(x) = \frac{\sum_{i=1}^{n} Y_i K((X_i-x)/h_n)}{\sum_{i=1}^{n} K((X_i-x)/h_n)},\] 
with some kernel $K(\cdot)$ and bandwidth $h_n\searrow 0$ (for $n \rightarrow \infty$),
which was introduced by \citet{Nadaraya1964} and \citet{Watson1964} as a nonparametric regression estimator  in a model $Y_i=m(X_i)+\varepsilon_i$ where the observations $(X_1,Y_1),\ldots,(X_n,Y_n)$ are assumed to be continuous. This estimator can be adapted for the classification case with categorical response $Y$ to estimate the posterior probability $P_g(x)=P(Y=g|x)$ as
\[\hat{P}_g(x) = \frac{\sum_{i=1}^{n} I\{Y_i=g\} K((X_i-x)/h_n)}{\sum_{i=1}^{n} K((X_i-x)/h_n)},\] 
see for instance \citet{HasTibFri2009}. The method introduced in \citet{SelkGertheiss2022} extends these estimators to handle multiple functional and categorical predictors, see Section \ref{model}. In the paper at hand, uniform convergence rates are given for the regression and classification estimators. In addition to estimating the regression function, variable selection is of interest, i.\,e.\ separating the relevant predictors from the noise variables. This is done with weights (counterpart of the bandwidth), which are data-driven for each covariate. The size of the weights then indicates the relevance of the corresponding covariate. In the paper at hand it is shown that the weights for noise variables actually vanish asymptotically, while those for relevant predictors tend to infinity. 

Uniform convergence rates in similar settings can be found, for instance, in \citet{LiOuyang2005}, \citet{FerratyEtal2010} or \citet{BouzebdaNemouchi2020}. A model with multiple continuous and categorical covariates is considered in \citet{LiOuyang2005}. The authors prove a uniform rate of convergence for a Nadaraya-Watson type regression estimator similar to the one considered in the paper at hand. Functional covariates are considered in \citet{FerratyEtal2010}. For a regression model with a single functional covariate, uniform convergence of a kernel regression estimator is shown. A similar result is given in \citet{BouzebdaNemouchi2020} who deal with conditional U-statistics.

Asymptotic results for weights / bandwidths determined by cross-validation in similar settings can for instance be found in \citet{HallLiRacine2007} or \citet{OuyangLiRacine2009}.
 \citet{HallLiRacine2007} prove for a regression model with continuous and categorical covariates that irrelevant predictors are smoothed out by cross-validated bandwidth. In \citet{OuyangLiRacine2009} a regression model with only categorical covariates is considered. This case has to be handled separately and cannot be treated as a special case of a model with mixed covariates, since the convergence rates for the smoothing parameters are different from the results in  \citet{HallLiRacine2007} and in the paper at hand.

The rest of the paper is organised as follows. In Section \ref{model}, the estimation method and the determination of the weights is explained. Further, assumptions needed for the theoretical results are stated. Section \ref{asymp} gives the asymptotic results, details the consistency and uniform convergence rates for the regression and classification estimators, and states the convergence of the data-driven weights. In Section \ref{sim} some finite sample simulation results are presented. Section \ref{proofs} sums up the proofs and Section \ref{outlook} concludes with a short discussion.

%%%%%%%%%%%%%%%%%%%%%%%%%%%%%%%%%%%%%%%%%%%%%%%%%

\section{Model}\label{model}

Let $(\bX_1,Y_1),\ldots,(\bX_n,Y_n)$ be iid observations, with $p$-dimensional covariates \\$\bX_i=(X_{i1},\ldots,X_{ip})$ that contain functional and categorical values, and continuous or categorical scalar response $Y_i$. The aim is to predict $Y$ given some new observation $\bx$ by estimating $E[Y_i|\bX_i=\bx]$ in a nonparametric way.
In the case of a continuous response, the regression case 
$Y_i=m(\bX_i)+\varepsilon_i$, this means estimating the regression function $m(\bx)$ and in the case of a categorical response $Y_i\in\{1,\ldots,G\}$, the classification case,  estimating the posterior probability $P_g(\bx)=P(Y=g|\bx)$ for all $g\in\{1,\ldots,G\}$. The estimator introduced by \citet{SelkGertheiss2022} for this target is
\begin{equation*} \frac{\sum_{i=1}^{n} Y_i K(\omega_1d_1(X_{i1},x_1) + \ldots + \omega_pd_p(X_{ip},x_p))}{\sum_{i=1}^{n} K(\omega_1d_1(X_{i1},x_1) + \ldots + \omega_pd_p(X_{ip},x_p))},\end{equation*}
in the regression case, and the same with $Y_i$ replaced by $I\{Y_i=g\}$ in the classification case, for some kernel function $K$ and weights $\omega_j$ that are determined in a data-driven way. The distance measures $d_1(\cdot,\cdot),\ldots,d_p(\cdot,\cdot)$ can take on very different types depending on the corresponding covariates. An example with categorical predictors $X_{il},x_l \in \{1,\ldots,G_l\}$ would be
\begin{equation}\label{distcat}d_l(X_{il},x_l) = \left\{ \begin{array}{ll}
0 & \mbox{ if } X_{il} = x_l,\\
1 & \mbox{ if } X_{il} \neq x_l.
\end{array}
\right.\end{equation}
For functional $X_{ij},x_j \in L^2$ one can use for example
\begin{equation}\label{metric}
d_j(X_{ij},x_j) = \sqrt{\int_{\mathcal{D}_j} (X_{ij}(t) - x_j(t))^2 dt},
\end{equation}
where $\mathcal{D}_j$ is the domain of the functions $X_{ij},x_j$. 
\begin{rem}
The distance measures in \eqref{distcat} and \eqref{metric} are just examples for possible (semi-)\\metrics that can be applied. If the categorical covariates are ordinal e.\,g.\ the metric $d_l(X_{il},x_l) =|X_{il}-x_l|$ would be appropriate. For functional variables besides metrics like the $L^2$ metric in \eqref{metric} it is also possible to use semi-metrics, see \citet[chap.\ 3]{FerratyVieu} or \citet[sec.\ 2.3] {SelkGertheiss2022}.
\end{rem}

In the paper at hand it is assumed that $X_{i1},\ldots,X_{ip_\f}$ are functional observations, $X_{i(p_\f+1)},\ldots,X_{i(p_\f+p_\c)}$ are categorical observations where $p=p_\f+p_\c$, and
\beq &&\exists\ 0<q_{\text{fun}}\leq p_{\text{fun}},\ 0\leq q_{\text{cat}}\leq p_{\text{cat}} : \\
&&X_{i1},\ldots,X_{iq_{\text{fun}}},X_{i(p_\text{fun}+1)},\ldots,X_{i(p_\text{fun}+q_{\text{cat}})}\ \text{are relevant predictors}\\
&& X_{i(q_{\text{fun}}+1)},\ldots,X_{ip_{\text{fun}}},X_{i(p_{\text{fun}}+q_{\text{cat}}+1)},\ldots,X_{i(p_{\text{fun}}+p_{\text{cat}})}\  \text{ are noise.}\eeq
Thus both types of covariates may contain relevant predictors and noise but there has to be at least one relevant functional predictor.
Further, it is assumed that the kernel has a product form $$K_\bo(\bx,\bx'):=\prod_{j=1}^{p_\text{fun}}K_j(\omega_jd_j(x_j,x'_j))\cdot\prod_{j=p_\text{fun}+1}^{p_\text{fun}+p_\text{cat}}(k_j^{-\omega_j})^{d_j(x_j,x'_j)}$$ for kernel functions $K_j$ and $k_j>1$, see Assumption (A3).  This is a typical form of a multivariate Nadaraya-Watson estimator, compare e.\,g.\ \citet{HaerdleMueller1997} for continuous covariates, the seminal paper of \citet{AitchisonAitken1976} for categorical ones, or \citet{LiOuyang2005} for mixed continuous and categorical covariates. Note further that $K(\omega_1d_1(X_{i1},x_1) + \ldots + \omega_pd_p(X_{ip},x_p))=\prod_{j=1}^{p_\text{fun}}K(\omega_jd_j(x_j,x'_j))\cdot\prod_{j=p_\text{fun}+1}^{p_\text{fun}+p_\text{cat}}(e^{-\omega_j})^{d_j(x_j,x'_j)}$ if $K$ is the Picard-kernel $K(u)=\exp(-u)I\{u\geq 0\}$. 

The performance of the estimators
\begin{equation}\label{hatm} \hat{m}(\bx): = \frac{\sum_{i=1}^n Y_i K_{\bo}(\bX_i,\bx)}{\sum_{i=1}^n K_{\bo}(\bX_i,\bx)}\end{equation}
and
\begin{equation}\label{hatP} \hat P_{g}(\bx):=\frac{\sum_{i=1}^nI\{Y_i=g\}K_{\bo}(\bX_i,\bx)}{\sum_{ i=1}^nK_{\bo}(\bX_i,\bx)}\end{equation}
strongly depends on the weights $\omega_1,\ldots,\omega_p$, see Theorem \ref{consistency}. A popular way to determine those weights, or their counterpart the bandwidths, is by cross-validation. With
\[\hat{m}_{(-i)}(\bX_i) = \frac{\sum_{s\neq i} Y_s K_{\bo}(\bX_s,\bX_i)}{\sum_{s\neq i} K_{\bo}(\bX_s,\bX_i)},\]
\[\hat P_{g(-i)}(\bX_i)=\frac{\sum_{s\neq i}I\{Y_s=g\}K_{\bo}(\bX_s,\bX_i)}{\sum_{s\neq i}K_{\bo}(\bX_s,\bX_i)}\]
being the leave-one-out estimates one may estimate $\omega_1,\ldots,\omega_p$ by minimising
\[Q(\omega_1,\ldots,\omega_p)=\frac 1n\sum_{i=1}^n (Y_i-\hat m_{(-i)}(\bX_i))^2v(\bX_i),\]
in the regression case and
\[Q(\omega_1,\ldots,\omega_p)=\frac 1n\sum_{i=1}^n\sum_{g=1}^G (I\{Y_i=g\}-\hat P_{g(-i)}(\bX_i))^2v(\bX_i)\]
in the classification case, respectively, where $v(\bx)=I\{\bx\in\cS\}$ for some suitable subset $\cS$, see Remark \ref{entropy} and Assumption (A7). The minimising weights are denoted by $\hat\omega_1,\ldots,\hat\omega_p$. Corollary \ref{convweights} states that with this data-driven weights the noise variables vanish asymptotically from the estimators defined in \eqref{hatm} and \eqref{hatP}. On the one hand this is interesting in the sense of variable selection. On the other hand the simulations in \citet{SelkGertheiss2022} show that downgrading the irrelevant covariates leads to an improvement of the estimation performance, see also Section \ref{sim}.

Before stating the assumptions needed to prove Theorem \ref{consistency} and \ref{rate} let's introduce some notation:
\begin{itemize}
\item Let $\bbX_i:=(X_{i1},\ldots,X_{iq_{\text{fun}}},X_{i(p_\text{fun}+1)},\ldots,X_{i(p_\text{fun}+q_{\text{cat}})})$ denote the relevant predictors and  $\tilde\bX_i =\bX_i\backslash \bbX_i$ the noise variables.
\item Let $\bX_\f$ denote the functional elements of $\bX$ and $\bX_\c$ the categorical ones. All functional covariates take values in semi-metric spaces $\mathcal{F}_j$ and each categorical covariate $X_{ij}$ takes values in $\{1,\ldots,\tilde G_j\}$. For the functional covariates suitable subsets $\cS_j$ of $\mathcal{F}_j$ are considered, see Remark  \ref{entropy} and Assumption (A7).
\item The indices are combined as follows \\
$ J_\f=\{1,\ldots,p_\text{fun}\}$, $J_\c=\{p_\text{fun}+1,\ldots,p_\text{fun}+p_{\text{cat}}\}$,\\
 $\bar J_\f=\{1,\ldots,q_{\text{fun}}\}$, $\bar J_\c=\{p_\text{fun}+1,\ldots,p_\text{fun}+q_{\text{cat}}\}$,
$\bar J=\bar J_\f \cup \bar J_\c$,\\
$\tilde J_\f=\{q_{\text{fun}}+1,\ldots,p_\text{fun}\}$, $\tilde J_\c=\{p_\text{fun}+q_{\text{cat}}+1,\ldots,p_\text{fun}+p_{\text{cat}}\}$, $\tilde J=\tilde J_\f \cup \tilde J_\c$.
\item For the distance measures set
$d_{\text{fun}}(\bx,\bx')=\max_{j\in J_\f}(d_j(x_j,x'_j))$ and
  $d_{\text{cat}}(\bx,\bx')=\max_{j\in J_\c}(d_j(x_j,x'_j))$. The categorical distance measures are standardized such that $d_j(x_j,x_j')\geq 1$ if $d_j(x_j,x_j')\neq 0$.
\end{itemize}
Throughout the paper $C$ and $C'$ denote positive constants whose values can change from line to line.

\subsubsection*{Assumptions:}
 \begin{description}
 \item[(A0)] As already stated, there must be at least one functional relevant covariate, i.\,e.\ $q_\f>0$. Further the categorical covariates are distributed such that \\$P(\max_{j\in J\c}d_j(X_j,x_j)=0|\bX_\f)\geq C_\c$ for some $C_\c>0$, almost surely, for all $(x_j)_{j\in J_\c}\in \bigtimes_{j\in J_\c}\{1,\ldots,\tilde G_j\}$. \\
 In the regression case $E[Y|\bbX]=:m(\bbX)$ and it is assumed that the second moment of $\e:=Y-m(\bbX)$ is bounded. 

 \end{description}
 \begin{description}
 \item[(A1)] $(\bbX,Y)$ is independent of $\tilde\bX$.
 \end{description}
  \begin{description}
 \item[(A2)] There exists some $\beta>0$ such that the regression function is Lipschitz of order $\beta$
 \[|m(\bx)-m(\bx')|\leq L(d_{\text{fun}}(\bx,\bx')^{\beta}+d_{\text{cat}}(\bx,\bx'))\]
  with some constant $L<\infty$.
The same holds for the posterior probability $P_g(\bx)=P(Y=g|\bX=\bx)$ for all $g\in\{1,\ldots,G\}$.
 \end{description}
 \begin{description}
 \item[(A3)] The kernel is defined as $K_\bo(\bx,\bx')=\prod_{j=1}^{p_\text{fun}}K_j(\omega_jd_j(x_j,x'_j))\cdot\prod_{j=p_\text{fun}+1}^{p_\text{fun}+p_\text{cat}}(k_j^{-\omega_j})^{d_j(x_j,x'_j)}$ where the kernel functions $K_j$ have bounded support $[0,1]$ and $k_j>1$. Further, there exist constants $0<C,C'<\infty$ such that $C\leq K_j(u)\leq C'$ for all $u\in[0,1]$. Also the following Lipschitz condition applies: $|\prod_{j=1}^{p_\f}K_j(u_j)-\prod_{j=1}^{p_\f}K_j(u'_j)|\leq L_K\max_{j\in J_\f} |u_j-u'_j|$ for some $L_K<\infty$ and all $u_j,u_j'\in [0,1]$. Additionally it is assumed that all $K_j$ have a unique maximum in zero.
 
 \end{description}
  \begin{description}
 \item[(A4)] There exist monotone bounded functions $\phi$ and $\varphi\leq 1$ and constants $0<C,C'<\infty$ such that for all $\epsilon>0$ and all $x_j\in\cS_j$, $j=1,\ldots,p_\f$, it holds
$C\phi(\epsilon)\leq P(\max_{j\in \bar J_\f}(d_j(X_j,x_j)) \leq\epsilon)\leq C'\phi(\epsilon)$  and $C\varphi(\epsilon)\leq P(\max_{j\in \tilde J_\f}(d_j(X_j,x_j)) \leq\epsilon)\leq C'\varphi(\epsilon)$. 
Further, for every functional covariate $X_j$ there exists a function $\phi_j$ with $C\phi_j(\epsilon)\leq P(d_j(X_j,x_j) \leq\epsilon)\leq C'\phi_j(\epsilon)$ that is Lipschitz.
 \end{description}
 Before proceeding with Assumptions (A5)-(A7) note the following remarks.
 \begin{rem}\label{smallball}
 The probability bounded by the functions $\phi, \varphi, \phi_j$ in Assumption (A4) is called the small ball probability and is known to play an important role for asymptotic results in nonparametric statistics for functional variables \citep{FerratyVieu}. It is easy to see that the choice of the distance measures $d_j$ directly influences the small ball probability. Chapter 13 in \citet{FerratyVieu} is dedicated to this issue. For a suitable choice  of distance measures it is possible to get $\phi_j(\epsilon)\sim \epsilon$, see Section 13.3.3 in the aforementioned book. Thus for the multivariate setting this means $\phi(\epsilon)\sim \epsilon^{q_\f}$ for a suitable dependence structure between the covariates.
% Theorem \ref{consistency} gives the convergence rate for the regression / classification estimator which depends on the small ball probability as well as the Lipschitz coefficient $\beta$ that also depends on the choice of the distance measure. 
 \end{rem}
\begin{rem}\label{abstand}
$K_\bo(\bX_s,\bX_i)$ is zero if $d_\f(\bbX_s,\bbX_i)= \max_{j\in\bar J_\f}(d_j(X_{sj},X_{ij})) >1/\min_{j\in \bar J_\f}(\omega_j)$. 
 \end{rem}
 \begin{rem}\label{EWkern}
 For $j=1,\ldots,p_\f$ let $x_j\in\cS_j$, and for $j=p_\f+1,\ldots,p_\f+p_\c$ let $x_j\in\{1,\ldots\tilde G_j\}$. Then uniformly in $\bx=(x_1,\ldots,x_p)$ it holds for some constants $C,C'<\infty$
 \beq
 E[K_\bo(\bX,\bx)]&=& E\Big[\prod_{j=1}^{p_\f}K_j(\omega_j d_j(X_j,x_j))\cdot\prod_{j=p_\f+1}^{p_\f+p_\c}(k_j^{-\omega_j})^{d_j(X_j,x_j)}\Big] \\
 &\leq& E\Big[\prod_{j=1}^{p_\f}K_j(\omega_j d_j(X_j,x_j))\Big] \\
 &\leq& C E[I\{\omega_j d_j(X_j,x_j) \leq 1 \forall j\in J_\f\}]\eeq 
 and 
 \beq
 E[I\{\omega_j d_j(X_j,x_j)\leq 1 \forall j\in J_\f\}]&=&P(\omega_j d_j(X_j,x_j)\leq 1 \forall j\in J_\f)\\
 &=&P(\omega_j d_j(X_j,x_j)\leq 1 \forall j\in\bar J_\f)\cdot P(\omega_j d_j(X_j,x_j)\leq 1 \forall j\in\tilde J_\f)\\
 &=&P(\max_{j\in\bar J_\f}(\omega_jd_j(X_j,x_j))\leq 1)\cdot P(\max_{j\in\tilde J_\f}(\omega_jd_j(X_j,x_j))\leq 1)\\
 &\leq&P(\min_{j\in\bar J_\f}(\omega_j)\max_{j\in\bar J_\f}(d_j(X_j,x_j))\leq 1)\cdot 1\\
 &\leq& C'\phi((\min_{j\in\bar J_\f}(\omega_j))^{-1}).
 \eeq
Further, for some constants $C,C'>0$ it holds
\beq
E[K_\bo(\bX,\bx)]&=& E\Big[\prod_{j=1}^{p_\f}K_j(\omega_j d_j(X_j,x_j))\cdot\prod_{j=p_\f+1}^{p_\f+p_\c}(k_j^{-\omega_j})^{d_j(X_j,x_j)}\Big] \\
&=& E\Big[\prod_{j=1}^{p_\f}K_j(\omega_j d_j(X_j,x_j))\cdot E\Big[\prod_{j=p_\f+1}^{p_\f+p_\c}(k_j^{-\omega_j})^{d_j(X_j,x_j)}|\bX_\f\Big]\Big] \\
&\geq& C E[I\{\omega_j d_j(X_j,x_j)\leq 1 \forall j\in J_\f\}] C_\c
\eeq
and
  \beq
&& E[I\{\omega_j d_j(X_j,x_j)\leq 1 \forall j\in J_\f\}]\\
&=&P(\omega_j d_j(X_j,x_j)\leq 1 \forall j\in J_\f)\\
 &=&P(\max_{j\in \bar J_\f}(\omega_jd_j(X_j,x_j))\leq 1)\cdot P(\max_{j\in\tilde J_\f}(\omega_jd_j(X_j,x_j))\leq 1)\\
 &\geq&P(\max_{j\in \bar J_\f}(\omega_j)\max_{j\in \bar J_\f}(d_j(X_j,x_j))\leq 1)\cdot P(\max_{j\in\tilde J_\f}(\omega_j)\max_{j\in\tilde J_\f}(d_j(X_j,x_j))\leq 1)\\
 &\geq& C\phi((\max_{j\in \bar J_\f}(\omega_j))^{-1})\varphi((\max_{j\in \tilde J_\f}(\omega_j))^{-1}).
 \eeq
For the sake of readability set $\phi_1:=\phi((\min_{j\in\bar J_\f}(\omega_j))^{-1})$ and \\$\phi_2:= \phi((\max_{j\in \bar J_\f}(\omega_j))^{-1})\varphi((\max_{j\in \tilde J_\f}(\omega_j))^{-1})$ and note that
$\phi_2\leq\phi_1$. 
 \end{rem}
 Now let's state two more assumptions (A5), (A6), followed by a further remark and the last assumption (A7).
   \begin{description}
 \item[(A5)] The weights $\omega_j$ may not go too fast to zero or to infinity. In detail, the following has to be fulfilled: $\tilde\phi:=\frac{\phi_2^2}{\phi_1}$ tends to zero with growing sample size $n$ - this holds if $\phi((\max_{j\in \bar J_\f}\omega_j)^{-1})$ tends to zero - but $\tilde\phi n^{\frac 13}\log(n)^{-1}$ tends to infinity.
%\\$\tilde\phi^2=o(1)$, $(\log(n))^3/(n\tilde\phi^3)=o(1)$
% \\$\log(n)/n\phi_1^{-1}\leq 1$: $\frac{\log(n)}{n\phi_1}\leq \frac{\log(n)}{n\tilde\phi}\leq \frac{\log(n)^3}{n\tilde\phi^3}=o(1)$
% \\$(n\phi_1)^{-1}=o(1)$ siehe oben
% \\$\log(n)/n\phi_2^{-1}=O(1)$ siehe oben
% \\$\sqrt{\frac{\log(n)}{n\tilde\phi}}\phi_2^{-1}=o(1)$: $\sqrt{\frac{\log(n)}{n\tilde\phi}}\phi_2^{-1}\leq \sqrt{\frac{\log(n)}{n\tilde\phi}}\tilde\phi^{-1}=\sqrt{\frac{\log(n)}{n\tilde\phi^3}}=o(1)$
% 
 Further $\max_{j\in J_\f}\omega_j=O(\sqrt{n/\log(n)})$.
% \\$\log(n)/n\max_{j\in J_\f}\omega_j\leq 1$ for $n$ large enough
% \\$\max_{j\in J_\f}\omega_j\phi_1=O\Big(\sqrt {n\tilde\phi/\log(n)}\Big)$ $\Leftrightarrow$ $\max_{j\in J_\f}\omega_j=O\Big(\sqrt {n/\log(n)}\sqrt{\frac{\phi_2^2}{\phi_1^2}}\Big)=O(\sqrt{n/\log(n)})$
 The term $\zeta(\bbo):=((\min_{j\in\bar J_\f}\omega_j)^{-\beta}+\max_{j\in\bar J_\c}(k_j^{-\omega_j}))^2$ has to fulfill $\zeta(\bbo)^{-1}=o(n)$ and $\zeta(\bbo)=o(\sqrt{n\tilde\phi\phi_2^2/\log(n)})$.
% \\$n^{-\frac 12}\sqrt{\zeta(\bbo)}=o(\zeta(\bbo))$  $\Leftrightarrow$ $\zeta(\bbo)^{-\frac 12}=o(n^{\frac 12})$
% \\ $\sqrt{\frac{\log(n)}{n\tilde\phi}}\phi_2^{-1}\zeta(\bbo)=o(1)$ 
%\\ $n^{-\frac 12}\sqrt{\zeta(\bbo)}=o(1)$: $\zeta(\bbo)=o(\sqrt{n\tilde\phi\phi_2^2/\log(n)})=o(n^{\frac 12})$
% \\$n^{-1}\tilde\phi^{-1}\zeta(\bbo)=o(1)$: $\zeta(\bbo)=o(\sqrt{n\tilde\phi\phi_2^2/\log(n)})=o(n^{\frac 12}\tilde\phi^{\frac 12})=o(n\tilde\phi)$
  \end{description}
 \begin{description}
 \item[(A6)] For the moments of the response variable $Y$ it holds $\forall l\geq 1$ $\exists C_l<\infty$ such that $E[|Y|^l] \leq C_l$.
 \end{description}
 \begin{rem}\label{entropy}
Let $\mathcal F$ be some semi-metric space and $\cS$ some subset. A finite set of points $x_1,\ldots,x_N$ in $\mathcal{F}$ is called an $\epsilon$-net for $\cS$ with some $\epsilon>0$ if $\cS$ is covered by the open balls with radius $\epsilon$ and center points $x_1,\ldots,x_N$. The minimum number of open balls with radius $\epsilon$ needed to cover $\cS$ is called the covering number of $\cS$ and denoted by $N_{\cS}(\epsilon)$. \\ 
In the setting at hand, for every functional covariate, a subset $\cS_j$ of $\mathcal{F}_j$ is considered and covered by an $\epsilon$-net $\mathcal{D}^{\epsilon}_j$ with $|\mathcal{D}^{\epsilon}_j|=N_{\cS_j}(\epsilon)$. These nets are combined to $\mathcal D_{\epsilon}=(\bigtimes_{j\in J_\f}\mathcal D^{\epsilon}_j)\times(\bigtimes_{j\in J_\c}\{1,\ldots,\tilde G_j\})$, see the proof of Lemma \ref{randden}. The combined covering number $N(\mathcal D_{\epsilon})$ is then a product of the univariate covering numbers $N(\mathcal D_{\epsilon})=\prod_{j\in J_\f}N_{\cS_j}(\epsilon)\prod_{j\in J_\c}\tilde G_j$. The Kolmogorov entropy $\psi(\mathcal D_{\epsilon}):=\log(N(\mathcal D_{\epsilon}))$ is thus a sum of the univariate entropies.
 In \citet{FerratyEtal2010} several examples are given on the size of the entropy for subsets of functional spaces. It is clear that the size of the entropy strongly depends on the choice of the semi-metric -  the distance measures $d_j$ in the setting at hand. Example 4 in \citet{FerratyEtal2010} is thus especially interesting  since it allows to construct, in any case, a semi-metric with reasonably ``small'' entropy.
 \end{rem}
  \begin{description}
\item[(A7)] The Kolmogorov entropy of $\cS:=(\bigtimes_{j\in J_\f}\cS_j)\times(\bigtimes_{j\in J_\c}\{1,\ldots,\tilde G_j\})$ fulfills \\$\psi\left(\mathcal D_{\frac{\log(n)}n}\right)=O(\log(n))$.
\end{description}
The assumptions are very similar to those stated in \citet{HallLiRacine2007} and \citet{FerratyEtal2010}, extended to the functional and / or multivariate case. (A2), (A3) and (A6) are typical for nonparametric settings and (A4) is a standard approach in dealing with functional data. (A7) is not necessarily needed and could be replaced by some weaker condition on the convergence rate of the entropy when Assumption (A5) is adapted accordingly. Nevertheless, it is stated like this for the sake of readability. Assumption (A1) might appear restrictive but it is needed for the proof technique used and stated in the same way in \citet{HallLiRacine2007}.

%%%%%%%%%%%%%%%%%%%%%%%%%%%%%%%%%%%%%%%%%%%%%%%%%%%

\section{Asymptotic results}\label{asymp}

Theorem \ref{consistency} states that the regression and classification estimators $\hat m$ and $\hat P_g$  with suitable weights $\bo$ are consistent and gives uniform convergence rates. 

\begin{theo}\label{consistency}
Under Assumptions (A0), (A2)--(A7) it holds
\[\sup_{\bx\in\cS}\left| \frac{\sum_{s=1}^n Y_s K_{\bo}(\bX_s,\bx)}{\sum_{s=1}^n K_{\bo}(\bX_s,\bx)}-m(\bbx)\right|=O\left((\min_{j\in\bar J_\f}\omega_j)^{-\beta}+\max_{j\in\bar J_\c}(k_j^{-\omega_j})\right)+O_P\left(\sqrt{\frac{\log(n)}{n\tilde\phi}}\right)\]
as well as for all $g=1,\ldots,G$
\beq
&&\sup_{\bx\in\cS}\left| \frac{\sum_{s=1}^n I\{Y_s=g\} K_{\bo}(\bX_s,\bx)}{\sum_{s=1}^n K_{\bo}(\bX_s,\bx)}-P_g(\bbx)\right|\\
&=&O\left((\min_{j\in\bar J_\f}\omega_j)^{-\beta}+\max_{j\in\bar J_\c}(k_j^{-\omega_j})\right)+O_P\left(\sqrt{\frac{\log(n)}{n\tilde\phi}}\right)\eeq
with $\tilde \phi=\frac{\phi_2^2}{\phi_1}$.
\end{theo}
The proof is given in Section \ref{proofs}. The result is comparable to those of \citet{LiOuyang2005} and \citet{FerratyEtal2010}, who give uniform convergence rates for a model with multiple continuous and categorical covariates and with a single functional covariate, respectively.

Theorem \ref{rate}  gives an asymptotic expansion  of  $Q$. With that the asymptotic behavior of $Q$`s minimising weights can be deduced, see Corollary \ref{convweights}.
\begin{theo}\label{rate}
Under Assumptions (A0)--(A7) with $\tilde \phi=\frac{\phi_2^2}{\phi_1}$ and $\zeta(\bbo)=((\min_{j\in\bar J_\f}\omega_j)^{-\beta}+\max_{j\in\bar J_\c}(k_j^{-\omega_j}))^2$ it holds
\beq
(i)\qquad Q(\bo)&=&A_m(\bbo)
+B_m(\bo)\\
&&+o_P(\zeta(\bbo)+(n\tilde\phi)^{-1}) +R
\eeq
and uniformly in $\bo$
\beq
(ii)\qquad Q(\bo)&=&A_m(\bbo)+o_P(1)+R
\eeq
 where $A_m(\bbo):=\frac{n-2}{n-1}E\Big[\Big(\frac{E\big[m(\bbX_1) \bar K_{\bbo}(\bbX_1,\bbX_2)|\bbX_2\big]}{E\big[ \bar K_{\bbo}(\bbX_1,\bbX_2)|\bbX_2\big]}-m(\bbX_2)\Big)^2v(\bX_2)\Big]=O(\zeta(\bbo))$, $B_m(\bo):=\frac {E[\e_1^2]}{n-1} E\Big[\frac{E\big[\big( K_{\bo}(\bX_1,\bX_2)\big)^2|\bX_2\big]}{\Big(E\big[ K_{\bo}(\bX_1,\bX_2)|\bX_2\big]\Big)^2}v(\bX_2)\Big]=O((n\tilde\phi)^{-1})$ and $R=O_P(1)$ is some term that does not depend on $\bo$. Accordingly in the classification case
\beq
(iii)\qquad Q(\bo)&=&A_P(\bbo)
+B_P(\bo)\\
&&+o_P(\zeta(\bbo)+(n\tilde\phi)^{-1}) +R'
\eeq
and uniformly in $\bo$
\beq
(iv)\qquad Q(\bo)&=&A_P(\bbo)+o_P(1)+R'
\eeq
where $A_P(\bbo):=\frac{n-2}{n-1}\sum_{g=1}^GE\Big[\Big(\frac{E\big[P_g(\bbX_1) \bar K_{\bbo}(\bbX_1,\bbX_2)|\bbX_2\big]}{E\big[ \bar K_{\bbo}(\bbX_1,\bbX_2)|\bbX_2\big]}-P_g(\bbX_2)\Big)^2v(\bX_2)\Big]=O(\zeta(\bbo))$, $B_P(\bo):=\frac {\sum_{g=1}^GE[(I\{Y_1=g\}-P_g(\bbX_1))^2]}{n-1} E\Big[\frac{E\big[\big( K_{\bo}(\bX_1,\bX_2)\big)^2|\bX_2\big]}{\Big(E\big[ K_{\bo}(\bX_1,\bX_2)|\bX_2\big]\Big)^2}v(\bX_2)\Big]=O((n\tilde\phi)^{-1})$ and $R'=O_P(1)$ is again some term that does not depend on $\bo$.
\end{theo}
The proof is given in Section \ref{proofs}. The result is comparable to that of \citet{HallLiRacine2007} who give a similar stochastic expansion for $Q$ in a model with multiple continuous and categorical covariates.

Corollary \ref{convweights} is a direct consequence of Theorem \ref{rate}. It states that the weights for relevant predictors converge to infinity whereas the weights for the noise variables vanish asymptotically.
\begin{cor}\label{convweights}
Assume that $A_m$ (resp.\ $A_P$) from Theorem \ref{rate} is a function of $(\omega_j)_{j\in\bar J}$ that goes to zero if and only if all $\omega_j$ for $j\in\bar J$ go to infinity. Then in both the regression and the classification case the minimiser $(\hat\omega_1,\ldots,\hat\omega_p)$ of $Q(\omega_1,\ldots,\omega_p)$ fulfills
\begin{itemize}
\item $\hat\omega_j\nPto\infty$ for all $j\in \bar J$ by Theorem \ref{rate} (ii) (resp.\ (iv)).
\item $\hat\omega_j\nPto 0$ for all $j\in\tilde J$ by Theorem \ref{rate} (i) (resp.\ (iii)) and Lemma \ref{eps2}.
\end{itemize}
Minimising $Q(\omega_1,\ldots,\omega_p)$ is a tradeoff between $\zeta(\bbo)$ and $(n\tilde\phi)^{-1}$ where the first goes faster to zero the faster $\bbo$ goes to infinity and the latter goes slower to zero the faster $\bbo$ goes to infinity. The rate of $(\hat\omega_j)_{j\in\bar J}$ as well as the rate in Theorem \ref{consistency} depends on the small ball probability $\phi$ which depends on the choice of the (semi-)metric $d_\f$, see Remark \ref{smallball}. With $\phi(\epsilon)\sim \epsilon^{q_\f}$ it holds 
\begin{itemize}
\item $\hat\omega_j\sim n^{\frac 1{q_\f+2\beta}}$ for all $j\in\bar J_\f$.
\item $\hat\omega_j\sim \frac \beta{q_\f+2\beta}\log(n)$ for all $j\in\bar J_\c$.
\item The convergence rate for the regression / classification estimator in Theorem \ref{consistency} with weights $(\hat\omega_1,\ldots,\hat\omega_p)$ is $O(n^{-\frac \beta{q_\f+2\beta}}\sqrt{\log(n)})$.
\end{itemize}
\end{cor}
The assumption on $A_m$ (resp.\ $A_P$) is stated in a very similar way in \citet{HallLiRacine2007}.

%%%%%%%%%%%%%%%%%%%%%%%%%%%%%%%%%%%%%%%%%%%%%%%

\section{Finite sample performance}\label{sim}

The performance of the estimators defined in \eqref{hatm} and \eqref{hatP} on simulated and real data is studied in detail in \citet{SelkGertheiss2022}. Here, only a small overview is given to illustrate the influence of a growing sample size and the downgrading of irrelevant covariates. Therefore, only a regression model with mixed covariates is presented and reference is made to \citet[sec.\ 3, 4]{SelkGertheiss2022} for further models.

For the regression model data is generated according to the following  model \citep[sec.\ 3]{SelkGertheiss2022}.
 For $i=1,\ldots,n$, $j=1\ldots,p_{\text{fun}}$, functional covariates are generated as
$$\tilde X_{ij}(t)=\sum_{l=1}^5\left(B_{ij,l}\sin\left(\frac tT(5-B_{ij,l})2\pi\right)-M_{ij,l}\right),$$ 
where $B_{ij,l}\sim\mathcal{U}[0,5]$ and $M_{ij,l}\sim\mathcal{U}[0,2\pi]$ for $l=1,\ldots,5$ and $T=300$.  $\mathcal{U}$ stands for the (continuous) uniform distribution.  $X_{1j}(t),\ldots,X_{nj}(t)$ are then determined as a centred and scaled version of $\tilde X_{1j}(t),\ldots,\tilde X_{nj}(t)$.
The categorical covariates are generated as
$X_{i(p_{\text{fun}}+1)},\ldots,X_{i(p_{\text{fun}}+p_{\text{cat}})}$ $\sim${\it B}$(0.5)$. With these functional and categorical covariates the regression is built as an extended functional linear model
$$Y_i=5\sum_{j=1}^{q_{\text{fun}}}\int X_{ij}(t)\gamma_{3,\frac 13}(t/10)dt+2(X_{i(p_{\text{fun}}+1)}+\ldots+X_{i(p_{\text{fun}}+q_{\text{cat}})})+\varepsilon_i$$
for some $q_{\text{fun}}\leq p_{\text{fun}}$ and $q_{\text{cat}}\leq p_{\text{cat}}$, where the coefficient function $\gamma_{a,b}(t)=b^a/\Gamma(a)t^{a-1}e^{-bt}I\{t>0\}$ is the density of the Gamma distribution. The errors $\varepsilon_i$ are iid standard normal. 
Two different cases are considered, the `minimal' and the `sparse' case, where $q_{\text{fun}}=q_{\text{cat}}=1$, $p_{\text{fun}}=p_{\text{cat}}=2$ in the minimal and $q_{\text{fun}}=q_{\text{cat}}=2$, $p_{\text{fun}}=p_{\text{cat}}=8$ in the sparse case. For all generated data sets a one-sided Picard kernel $K(u)=e^{-u}I\{u\geq 0\}$ is used and the results shown are based on 500 replications each.

Additionally to the method explained in Section \ref{model} the minimiser of $Q$ is determined under the restrictions
\begin{itemize}
\item[(i)] $\omega_1=\omega_2=\ldots=\omega_p$,
\item[(ii)] $\omega_j=0$ for all covariates with no influence on the response. %$j$ with $x_j$ has no influence on $y$.
\end{itemize}
Under restriction (ii), which is also referred to as the `oracle', the minimising weights are thus determined only for the relevant covariates. Restriction (i), on the other hand, leads to a single minimising weight and can be interpreted as determining an appropriate overall/global bandwidth. A comparison of the results under both restrictions also shows the difference in estimation performance with and without inclusion of the irrelevant covariates in the estimator.

\begin{figure}
\begin{centering}
\subfigure{\includegraphics[width=0.5\textwidth]{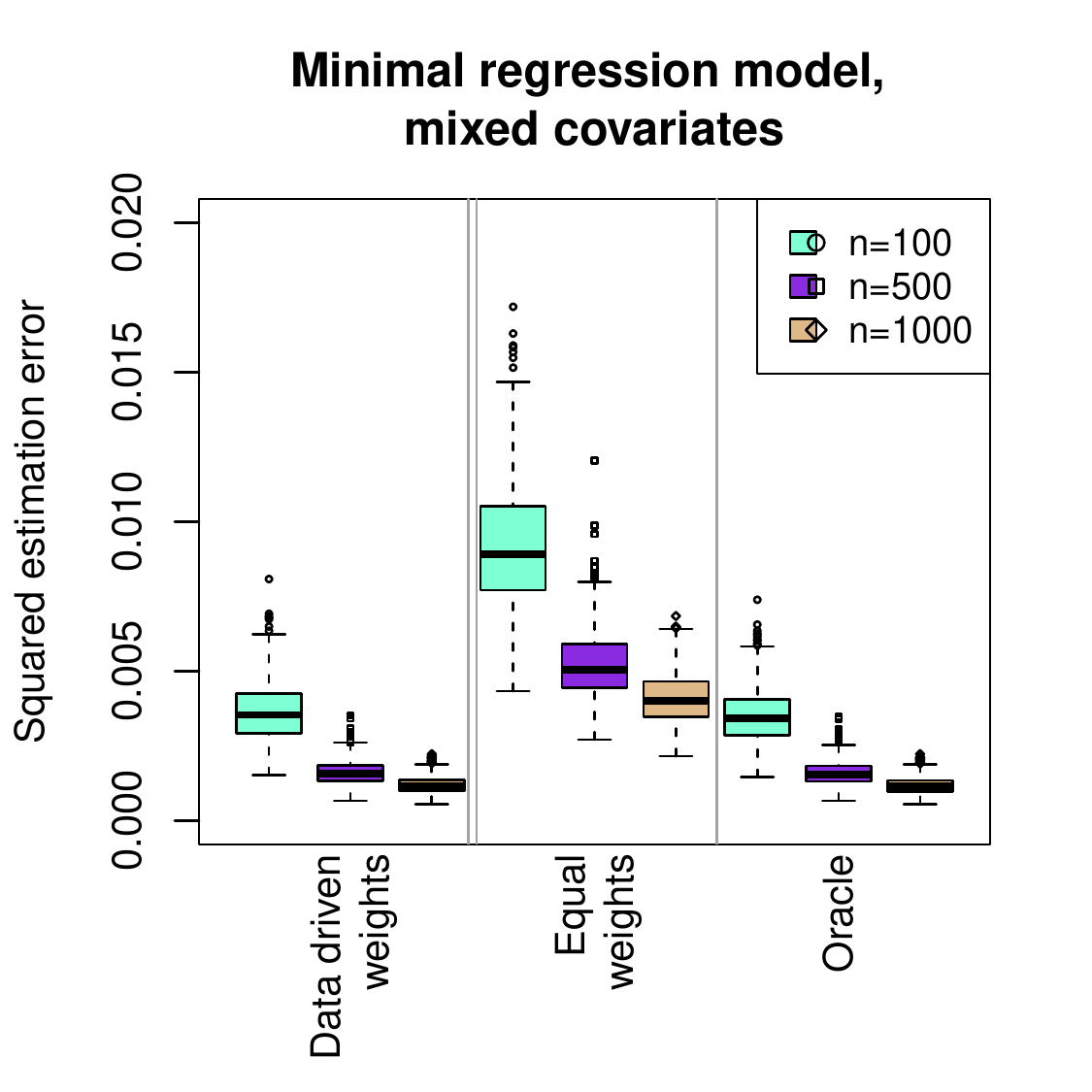}}\hspace{-0.02\textwidth}
\subfigure{\includegraphics[width=0.5\textwidth]{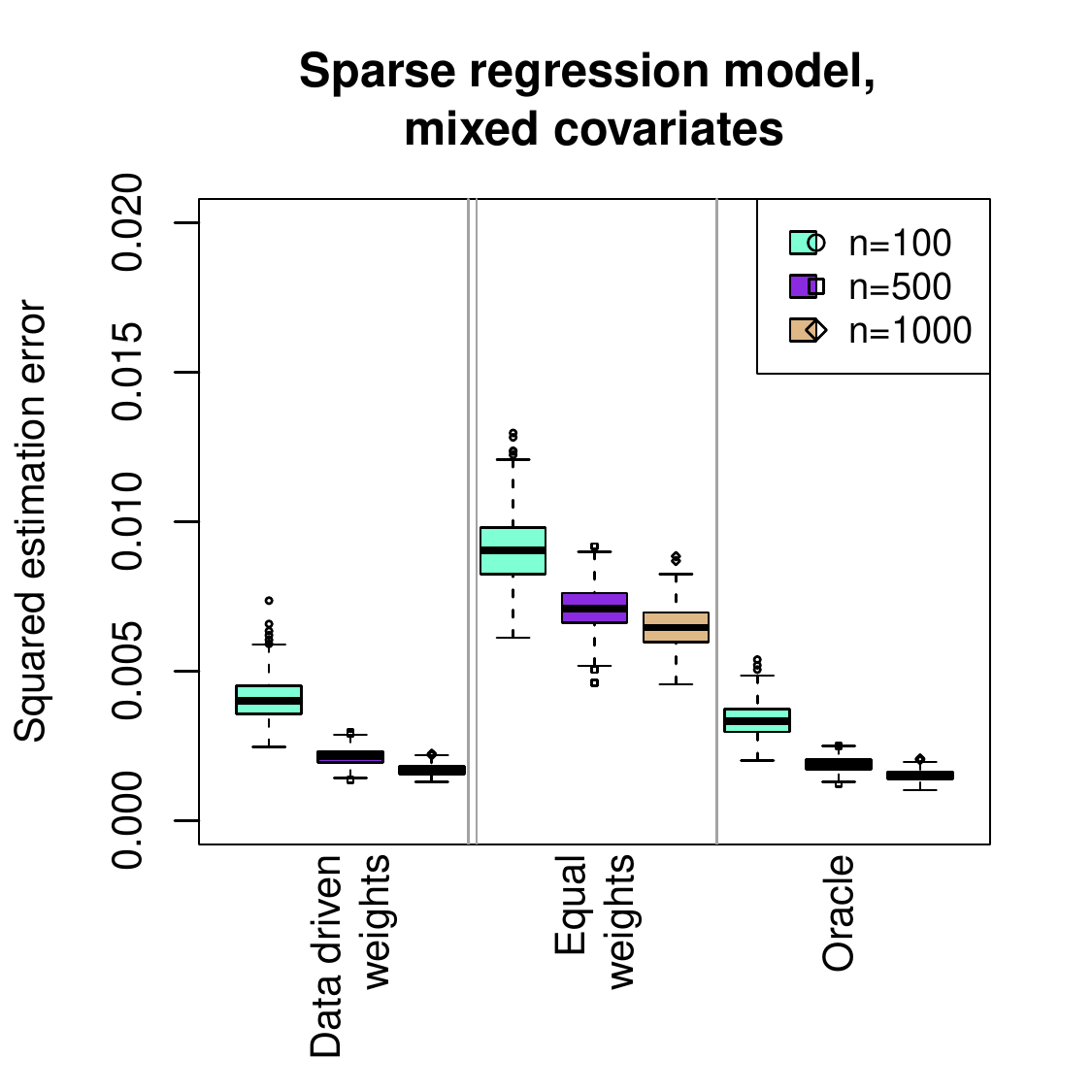}}
\end{centering}
\caption{Prediction results  in the minimal (left) and sparse (right) case with no restriction (`data driven weights'), restriction (i, `equal weights') and (ii, `oracle'), respectively. }\label{fig:Remse}
\end{figure}

\begin{figure}
\begin{centering}
\subfigure{\includegraphics[width=0.41\textwidth]{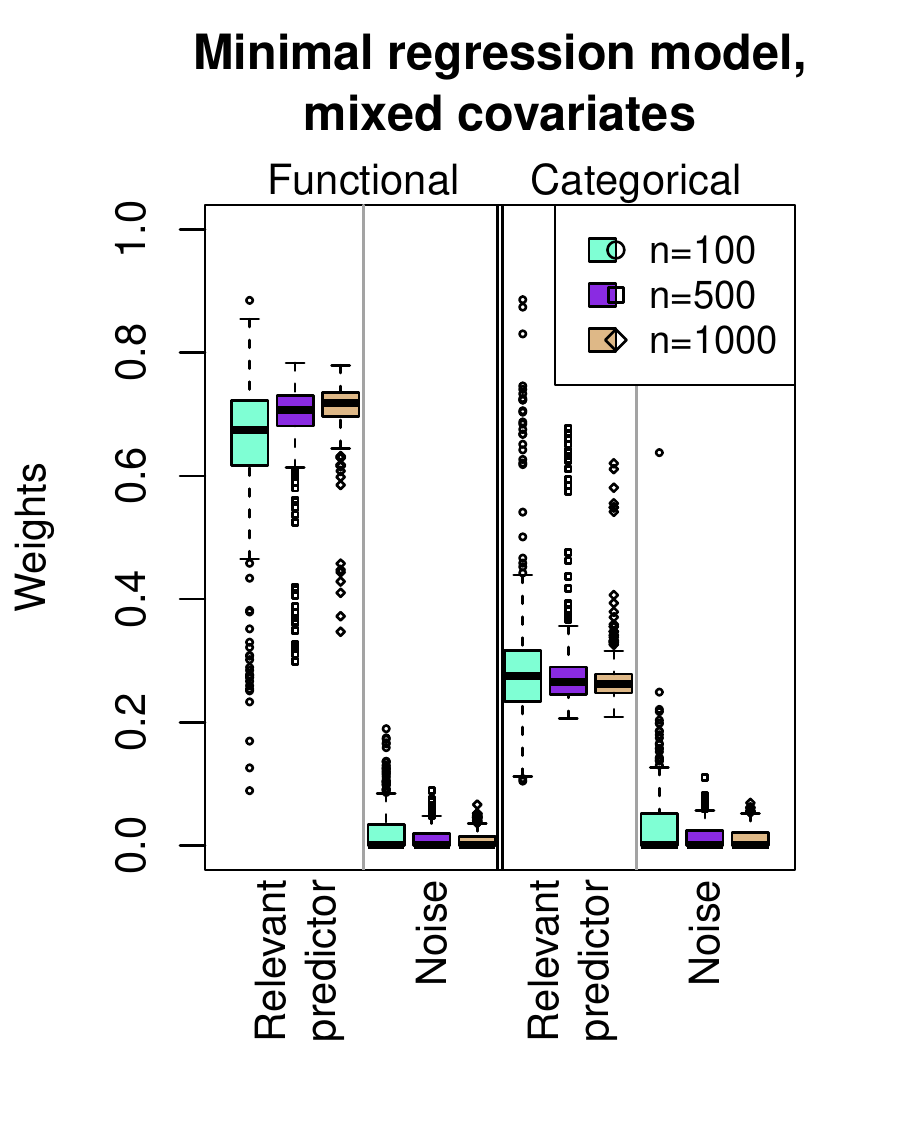}}\hspace{-0.02\textwidth}
\subfigure{\includegraphics[width=0.61\textwidth]{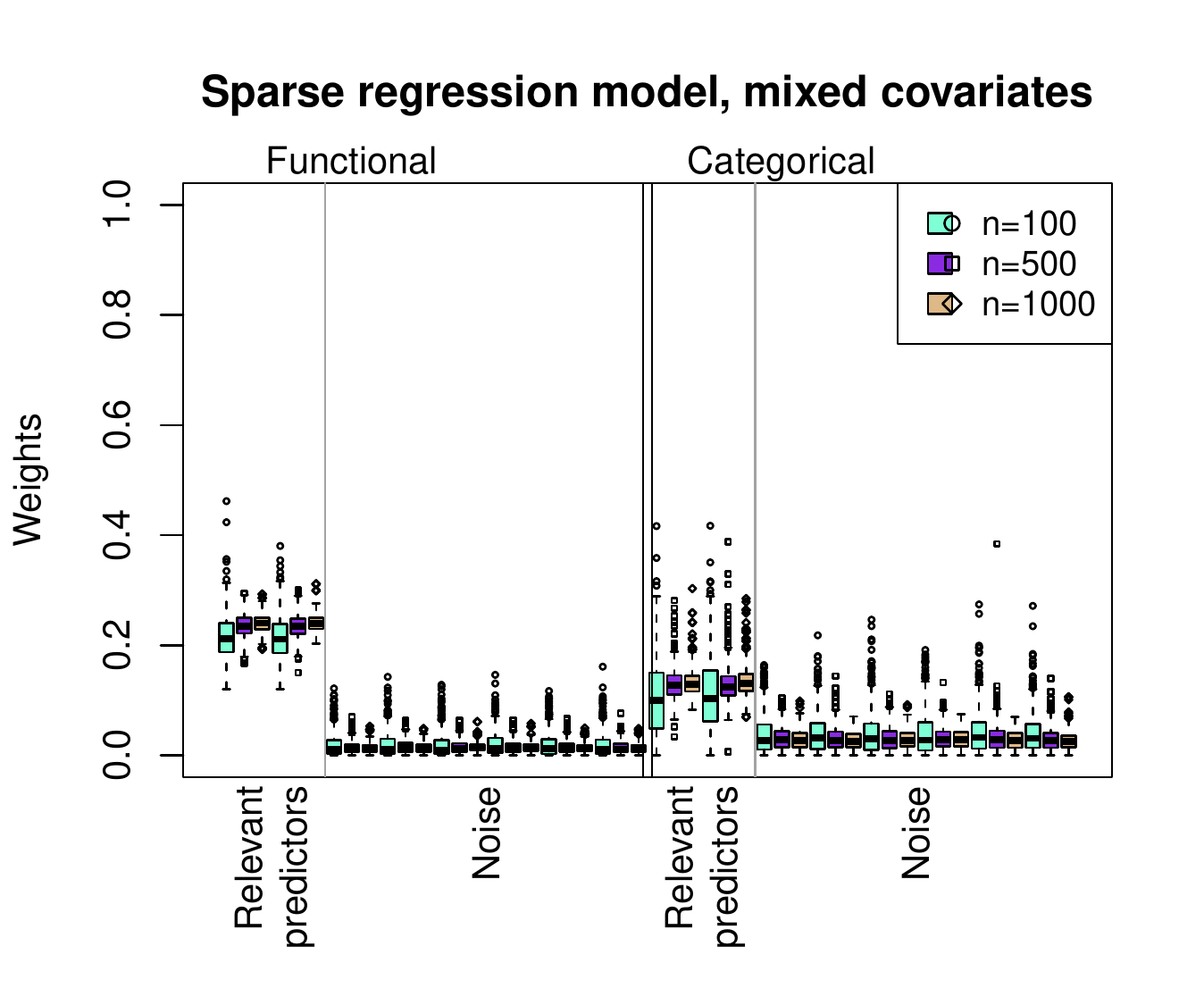}}
\end{centering}
\caption{Normed minimising weights $\frac{\hat\omega_j}{\sum_{k=1}^p\hat\omega_k}$ in the minimal and sparse case, respectively. }\label{fig:Reweights}
\end{figure}

Figure \ref{fig:Remse} displays the estimation results for different sample sizes $n=100, 500, 1000$. It can be seen that the performance increases with growing sample size and that the results are better when the irrelevant variables are excluded from the estimator, where the data-driven weights achieve results comparable to the oracle case. In Figure \ref{fig:Reweights} the estimated weights are presented for different sample sizes $n=100, 500, 1000$. The effect that (only) the weights for the noise variables vanish with growing sample size can clearly be seen.

In the classification case, the performance of the estimator \eqref{hatP} is also very good. In \citet[sec. 3.2.2]{SelkGertheiss2022} this is shown for different models using the misclassification rate and the squared estimation error.

%%%%%%%%%%%%%%%%%%%%%%%%%%%%%%%%%%%%%%%
%%%%%%%%%%%%%%%%%%%%%%%%%%%%%%%%%%%%%%%

\section{Proofs}\label{proofs}
{\bf Proof of Theorem \ref{consistency}:}
Consider the decomposition
\beq
 &&\frac{\sum_{s=1}^n Y_s K_{\bo}(\bX_s,\bx)}{\sum_{s=1}^n K_{\bo}(\bX_s,\bx)}-m(\bbx)\\
 &=&\frac {E[K_\bo(\bX_1,\bx)]}{\frac 1n\sum_{s=1}^n K_{\bo}(\bX_s,\bx)}\Big(\frac{\frac1{n}\sum_{s=1}^n Y_sK_\bo(\bX_s,\bx)}{E[K_\bo(\bX_1,\bx)]}-\frac{E\big[\frac1{n}\sum_{s=1}^n Y_sK_\bo(\bX_s,\bx)\big]}{E[K_\bo(\bX_1,\bx)]}\Big)\\
 &&+\frac {E[K_\bo(\bX_1,\bx)]}{\frac 1n\sum_{s=1}^n K_{\bo}(\bX_s,\bx)}\Big(\frac{E\big[\frac1{n}\sum_{s=1}^n Y_sK_\bo(\bX_s,\bx)\big]}{E[K_\bo(\bX_1,\bx)]}-m(\bbx)\Big)\\
 &&+\frac {E[K_\bo(\bX_1,\bx)]}{\frac 1n\sum_{s=1}^n K_{\bo}(\bX_s,\bx)}\Big(1-\frac{\frac1{n}\sum_{s=1}^n K_\bo(\bX_s,\bx)}{E[K_\bo(\bX_1,\bx)]}\Big)m(\bbx).
\eeq
For the classification case the same holds with $Y_s$ replaced by $I\{Y_s=g\}$ and $m(\bbx)$ replaced by $P_g(\bbx)$.
Now the assertion can be deduced from the following Lemmas \ref{randden}, \ref{Em} and \ref{EY}. Note that $\sup_{\bx\in\cS}|m(\bbx)|<\infty$ and that Lemma \ref{randden} implies $\frac {E[K_\bo(\bX_1,\bx)]}{\frac 1n\sum_{s=1}^n K_{\bo}(\bX_s,\bx)}=O_P(1)$ uniformly in $\bx\in\cS$.

%%%%%%%%%%%%%%%%%%%%%%%%%%%%%%%%%%%

\begin{lemma}\label{randden}
Under Assumptions (A3), (A4), (A5) and (A7) it holds
\[\sup_{\bx\in\cS}\Big|\frac{\frac1{n}\sum_{s=1}^n K_\bo(\bX_s,\bx)}{E[K_\bo(\bX_1,\bx)]}-1\Big|=O_P\left(\sqrt{\frac{\log(n)}{n\tilde\phi}}\right)\]
with $\tilde \phi=\frac{\phi_2^2}{\phi_1}$.
\end{lemma}
{\bf Proof of Lemma \ref{randden}:}
Following the proof of Lemma 8 in \citet{FerratyEtal2010} set $\epsilon=\frac{\log(n)}n$ and let $\mathcal D^{\epsilon}_j$, $j=1,\ldots,p_\f$ be minimal $\epsilon$-nets for $\mathcal S_j$, as defined in \citet{FerratyEtal2010} and Remark \ref{entropy}, and $\mathcal D_{\epsilon}=(\bigtimes_{j\in J_\f}\mathcal D^{\epsilon}_j)\times(\bigtimes_{j\in J_\c}\{1,\ldots,\tilde G_j\})$. For all $\bx$ set $\bx_{k(\bx)}=\arg\min_{\bx'\in\mathcal D_{\epsilon}}d(\bx,\bx')$. With $\psi(\mathcal D_{\epsilon})$ the Kolmogorov entropy is denoted and with $N(\mathcal D_{\epsilon})$ the covering number, see Remark \ref{entropy}. Then consider the decomposition
\[\sup_{\bx\in\cS}\Big|\frac{\frac1{n}\sum_{s=1}^n K_\bo(\bX_s,\bx)}{E[K_\bo(\bX_1,\bx)]}-1\Big|\leq\underbrace{\sup_{\bx\in\cS}|\hat f(\bx)-\hat f(\bx_{k(\bx)})|}_{F_1}+\underbrace{\sup_{\bx\in\cS}|\hat f(\bx_{k(\bx)})-1|}_{F_2}\]
where $\hat f(\bx):=\frac{\frac1{n}\sum_{s=1}^n K_\bo(\bX_s,\bx)}{E[K_\bo(\bX_1,\bx)]}$.
For some constant $C<\infty$ it holds with Remark \ref{EWkern}
\beq
F_1&\leq& \frac C{\phi_2}\sup_{\bx\in\cS}\frac 1n\sum_{s=1}^n|K_\bo(\bX_s,\bx)-K_\bo(\bX_s,\bx_{k(\bx)})|\\
&=&\frac C{\phi_2}\sup_{\bx\in\cS}\frac 1n\sum_{s=1}^n\Big|\prod_{j=1}^{p_\f}K_j(\omega_j d_j(X_{sj},x_j))-\prod_{j=1}^{p_\f}K_j(\omega_j d_j(X_{sj},x_{k(\bx)j}))\Big|\cdot\prod_{j=p_\f+1}^{p_\f+p_\c}(k_j^{-\omega_j})^{d_j(X_{sj},x_j)}\\
&\leq&\frac C{\phi_2}\sup_{\bx\in\cS}\frac 1n\sum_{s=1}^n\Big|\prod_{j=1}^{p_\f}K_j(\omega_j d_j(X_{sj},x_j))-\prod_{j=1}^{p_\f}K_j(\omega_j d_j(X_{sj},x_{k(\bx)j}))\Big|\\
&=&\frac C{\phi_2}\sup_{\bx\in\cS}\frac 1n\sum_{s=1}^n\Big|\prod_{j=1}^{p_\f}K_j(\omega_j d_j(X_{sj},x_j))-\prod_{j=1}^{p_\f}K_j(\omega_j d_j(X_{sj},x_{k(\bx)j}))\Big|\\
&&\hspace{4cm}\cdot I\{(\omega_j d_j(X_{sj},x_{j})\leq 1\forall j\in J_\f)\vee(\omega_j d_j(X_{sj},x_{k(\bx)j})\leq 1\forall j\in J_\f)\}.
\eeq
Now consider the decomposition
\[F_1\leq C'\sup_{\bx\in\cS}(F_{11}+F_{12}+F_{13})\]
with 
\beq
F_{11}&=&\frac 1{\phi_2}\sup_{\bx\in\cS}\frac 1n\sum_{s=1}^n\Big|\prod_{j=1}^{p_\f}K_j(\omega_j d_j(X_{sj},x_j))-\prod_{j=1}^{p_\f}K_j(\omega_j d_j(X_{sj},x_{k(\bx)j}))\Big|\\
&&\hspace{4cm}\cdot I\{(\omega_j d_j(X_{sj},x_{j})\leq 1\forall j\in J_\f)\wedge(\omega_j d_j(X_{sj},x_{k(\bx)j})\leq 1\forall j\in J_\f)\} \\\\
F_{12}&=&\frac 1{\phi_2}\sup_{\bx\in\cS}\frac 1n\sum_{s=1}^n I\{(\omega_j d_j(X_{sj},x_{j})\leq 1\forall j\in J_\f)\wedge(\exists j\in J_\f: \omega_j d_j(X_{sj},x_{k(\bx)j})> 1)\}\\\\
F_{13}&=&\frac 1{\phi_2}\sup_{\bx\in\cS}\frac 1n\sum_{s=1}^n I\{(\omega_j d_j(X_{sj},x_{k(\bx)j})\leq 1\forall j\in J_\f)\wedge(\exists j\in J_\f: \omega_j d_j(X_{sj},x_{j})> 1)\}.
\eeq
By the Lipschitz assumption on $K$ in (A3) it holds for some constants $C,C'<\infty$
\beq
F_{11}&\leq& \frac {C}{\phi_2}\sup_{\bx}\frac 1n\sum_{s=1}^n\max_{j\in J_\f}|\omega_j (d_j(X_{sj},x_j))- d_j(X_{sj},x_{k(\bx)j}))|\\
&&\hspace{4cm}\cdot I\{(\omega_j d_j(X_{sj},x_{j})\leq 1\forall j\in J_\f)\wedge(\omega_j d_j(X_{sj},x_{k(\bx)j})\leq 1\forall j\in J_\f)\} \\
&\leq&\frac {C'\epsilon\max_{j\in J_\f}\omega_j}{\phi_2}\sup_{\bx\in\cS}\frac 1n\sum_{s=1}^n\Big|I\{(\omega_j d_j(X_{sj},x_{j})\leq 1\forall j\in J_\f)\wedge(\omega_j d_j(X_{sj},x_{k(\bx)j})\leq 1\forall j\in J_\f)\}\Big|
\eeq
and thus with a Bernstein inequality  \citep[Corollary A.9]{FerratyVieu} for all $\eta>0$ and some constants $0<C,C'<\infty$
%$O(\frac {\epsilon\max_{j\in J_\f}\omega_j}{\phi_2})$, $E=O(\frac {\epsilon\max_{j\in J_\f}\omega_j\phi_1}{\phi_2})$, $\Var=O(\frac {\epsilon^2\max_{j\in J_\f}\omega_j^2\phi_1}{\phi_2^2})$
\beq
P\left(|F_{11}-E[F_{11}]|>\eta \sqrt{\frac{\log(n)}{n\tilde\phi}}\right)&\leq&2\exp\left(-\frac{n\eta^2\frac{\log(n)}{n\tilde\phi}}{2C\left(\frac {\epsilon^2\max_{j\in J_\f}\omega_j^2\phi_1}{\phi_2^2}+\sqrt{\frac{\log(n)}{n\tilde\phi}} \frac {\epsilon\max_{j\in J_\f}\omega_j}{\phi_2}\right)}\right)\\
&\leq&2\exp\big(-C'\eta^2\log(n)\big)\\
&=&o(1)
\eeq
since $\phi_2^{-1}\leq \tilde\phi^{-1}$ and $\epsilon\max_{j\in J_\f}\omega_j\leq 1$ for $n$ large enough. Combine this result with $E[F_{11}]=O(\frac {\epsilon\max_{j\in J_\f}\omega_j\phi_1}{\phi_2})$ to get
\[F_{11}=O\Big(\frac {\epsilon\max_{j\in J_\f}\omega_j\phi_1}{\phi_2}\Big)+O_P\left(\sqrt{\frac{\log(n)}{n\tilde\phi}}\right)=O_P\left(\sqrt{\frac{\log(n)}{n\tilde\phi}}\right)\]
with $\max_{j\in J_\f}\omega_j\phi_1=O\Big(\sqrt {n\tilde\phi/\log(n)}\Big)$.
For $F_{12}$ it holds 
\beq
&&E\big[I\{(\omega_j d_j(X_{sj},x_{j})\leq 1\forall j\in J_\f)\wedge(\exists j\in J_\f: \omega_j d_j(X_{sj},x_{k(\bx)j})> 1)\}\big]\\
&=&P\Big((\omega_j d_j(X_{sj},x_{j})\leq 1\forall j\in J_\f)\cap (\exists j\in J_\f: \omega_j d_j(X_{sj},x_{k(\bx)j})> 1)\Big)\\
&\leq&P\Big(\exists j\in J_\f: \omega_j d_j(X_{sj},x_{j})\leq 1\wedge  \omega_j d_j(X_{sj},x_{k(\bx)j})> 1\Big)\\
&\leq&\sum_{j\in J_\f}P\big(\omega_j d_j(X_{sj},x_{j})\leq 1\cap  \omega_j d_j(X_{sj},x_{k(\bx)j})> 1\big)\\
&\leq&\sum_{j\in J_\f}P\big(\omega_j d_j(X_{sj},x_{j})\leq 1\cap  \omega_j (d_j(X_{sj},x_{j})+\epsilon)> 1\big)\\
&=& \sum_{j\in J_\f}\Big(P\big(\omega_j d_j(X_{sj},x_{j})\leq 1\big)-P\big(\omega_j d_j(X_{sj},x_{j})\leq 1-\omega_j\epsilon\big)\Big)\\
&\leq&\sum_{j\in J_\f}C(\phi_j(\omega_j^{-1})-\phi_j(\omega_j^{-1}-\epsilon))\\
&=&O(\epsilon)
\eeq
and for all $\eta>0$ and some constants $0<C,C'<\infty$ 
%$O(\frac {1}{\phi_2})$, $E=O(\frac {\epsilon}{\phi_2})$, $\Var=O(\frac {\epsilon}{\phi_2^2})$
\beq
P\left(|F_{12}-E[F_{12}]|>\eta  \sqrt{\frac{\log(n)}{n\tilde\phi}}\right)&\leq&2\exp\left(-\frac{n\eta^2\frac{\log(n)}{n\tilde\phi}}{2C\left(\frac {\epsilon}{\phi_2^2}+\sqrt{\frac{\log(n)}{n\tilde\phi}} \frac {1}{\phi_2}\right)}\right)\\
&\leq&2\exp\big(-C'\eta^2\log(n)\big)\\
&=&o(1)
\eeq
with $\frac{\epsilon}{\phi_1}\leq 1$ for $n$ large enough. Thus
\[F_{12}=O_P\left(\sqrt{\frac{\log(n)}{n\tilde\phi}}\right).\] 
With the same arguments one gets $F_{13}=O_P\left(\sqrt{\frac{\log(n)}{n\tilde\phi}}\right)$.

Considering $F_2$ for all $\eta>0$ it holds
\beq
P\left(F_2>\eta\sqrt{\frac{\psi(\mathcal D_\epsilon)}{n\tilde\phi}}\right)&=&P\left(\max_{\bx'\in\mathcal D_{\epsilon}}|\hat f(\bx')-1|>\eta\sqrt{\frac{\psi(\mathcal D_\epsilon)}{n\tilde\phi}}\right)\\
&\leq&N(\mathcal D_{\epsilon})\max_{\bx'\in\mathcal D_{\epsilon}}P\left(|\hat f(\bx')-1|>\eta\sqrt{\frac{\psi(\mathcal D_\epsilon)}{n\tilde\phi}}\right).
\eeq
Further, $\frac{ K_\bo(\bX_s,\bx')}{E[K_\bo(\bX_1,\bx')]}=O(\phi_2^{-1})=O\left(\frac{\phi_1}{\phi_2^2}\right)$, $E\left[\frac{ K_\bo(\bX_s,\bx')}{E[K_\bo(\bX_1,\bx')]}\right]=1$ and $\Var\left(\frac{ K_\bo(\bX_s,\bx')}{E[K_\bo(\bX_1,\bx')]}\right)=O\left(\frac{\phi_1}{\phi_2^2}\right)+O(1)=O\left(\frac{\phi_1}{\phi_2^2}\right)$ for all $s=1,\ldots,n$ and $\bx'\in\mathcal D_\epsilon$ by Remark \ref{EWkern} and thus with a Bernstein inequality as before for some constants $0<C,C'<\infty$
\beq
P\left(|\hat f(\bx')-1|>\eta\sqrt{\frac{\psi(\mathcal D_\epsilon)}{n\tilde\phi}}\right)&=&P\left(\frac 1n\left|\sum_{s=1}^n\left(\frac{ K_\bo(\bX_s,\bx')}{E[K_\bo(\bX_1,\bx')]}-1\right)\right|>\eta\sqrt{\frac{\psi(\mathcal D_\epsilon)}{n\tilde\phi}}\right)\\
&\leq& 2\exp\left(-\frac{n\eta^2\frac{\psi(\mathcal D_\epsilon)}{n\tilde\phi}}{2 C(\frac{\phi_1}{\phi_2^2}+\sqrt{\frac{\psi(\mathcal D_\epsilon)}{n\tilde\phi}}\frac{\phi_1}{\phi_2^2})}\right)\\
&\leq& 2\exp\left(-C'\eta^2\psi(\mathcal D_\epsilon)\right).
\eeq
Further,
\beq
P\left(F_2>\eta\sqrt{\frac{\psi(\mathcal D_\epsilon)}{n\tilde\phi}}\right)&\leq&2N(\mathcal D_\epsilon)\exp\big(-C'\eta^2\psi(\mathcal D_\epsilon)\big)\\
&=&2\exp\big((1-C'\eta^2)\psi(\mathcal D_\epsilon)\big)
\eeq
which completes the proof since $\eta$ can be chosen such that $C'\eta^2>1$.

\boxi

%%%%%%%%%%%%%%%%%%%%%%%%%%%%%%%%%%%%%%%

\begin{lemma}\label{Em}
Under Assumptions (A0), (A2), (A3) and (A4) it holds
\[\sup_{\bx\in\cS}\Big|\frac{E\big[\frac1{n}\sum_{s=1}^n Y_sK_\bo(\bX_s,\bx)\big]}{E[K_\bo(\bX_1,\bx)]}-m(\bbx)\Big|=O\left((\min_{j\in\bar J_\f}\omega_j)^{-\beta}+\max_{j\in\bar J_\c}(k_j^{-\omega_j})\right)\]
as well as
\[\sup_{\bx\in\cS}\Big|\frac{E\big[\frac1{n}\sum_{s=1}^n I\{Y_s=g\}K_\bo(\bX_s,\bx)\big]}{E[K_\bo(\bX_1,\bx)]}-P_g(\bbx)\Big|=O\left((\min_{j\in\bar J_\f}\omega_j)^{-\beta}+\max_{j\in\bar J_\c}(k_j^{-\omega_j})\right)\]
for all $g=1,\ldots,G$.
\end{lemma}
{\bf Proof of Lemma \ref{Em}:}
It holds
\beq
\sup_{\bx\in\cS}\Big|\frac{E\big[\frac1{n}\sum_{s=1}^n Y_sK_\bo(\bX_s,\bbx)\big]}{E[K_\bo(\bX_1,\bx)]}-m(\bbx)\Big|&=&\sup_{\bx\in\cS}\Big|\frac{ E[K_\bo(\bX_1,\bx)(Y_1-m(\bbx))]}{E[K_\bo(\bX_1,\bx)]}\Big|\\
&=&\sup_{\bx\in\cS}\Big|\frac{ E[K_\bo(\bX_1,\bx)(m(\bbX_1)-m(\bbx))]}{E[K_\bo(\bX_1,\bx)]}\Big|\\
&\leq&\sup_{\bx\in\cS}\frac{ E[K_\bo(\bX_1,\bx)|m(\bbX_1)-m(\bbx)|]}{E[K_\bo(\bX_1,\bx)]},
\eeq
since the errors are centred and independent of $\bX$. For the classification case similarly
\beq
&&\sup_{\bx\in\cS}\Big|\frac{E\big[\frac1{n}\sum_{s=1}^n I\{Y_s=g\}K_\bo(\bX_s,\bbx)\big]}{E[K_\bo(\bX_1,\bx)]}-P_g(\bbx)\Big|\\
&=&\sup_{\bx\in\cS}\Big|\frac{ E[E[K_\bo(\bX_1,\bx)(I\{Y_1=g\}-P_g(\bbx))]|\bbX_1]}{E[K_\bo(\bX_1,\bx)]}\Big|\\
&=&\sup_{\bx\in\cS}\Big|\frac{ E[K_\bo(\bX_1,\bx)(P_g(\bbX_1)-P_g(\bbx))]}{E[K_\bo(\bX_1,\bx)]}\Big|\\
&\leq&\sup_{\bx\in\cS}\frac{ E[K_\bo(\bX_1,\bx)|P_g(\bbX_1)-P_g(\bbx)|]}{E[K_\bo(\bX_1,\bx)]}.
\eeq
Further, by the Lipschitz assumption on $m$ (same for $P_g$) and Remark \ref{abstand} one gets
\beq
\sup_{\bx\in\cS}\frac{ E[K_\bo(\bX_1,\bx)|m(\bbX_1)-m(\bx)|]}{E[K_\bo(\bX_1,\bx)]}&\leq&L\Big((\min_{j\in\bar J_\f}\omega_j)^{-\beta}+\sup_{\bx\in\cS}\frac{ E[K_\bo(\bX_1,\bx)d_\c(\bbX_1,\bbx)]}{E[K_\bo(\bX_1,\bx)]}\Big)\\
&=&O\Big((\min_{j\in\bar J_\f}\omega_j)^{-\beta}+\max_{j\in\bar J_\c}(k_j^{-\omega_j})\Big),
\eeq
where the last equality holds, since if $d_\c(\bbX_1,\bbx)\neq 0$ there exists some $j_1\in \bar J_\c$ with $d_{j_1}(X_{1j_1},x_{j_1}))\geq 1$ and thus  $k_{j_1}^{-\omega_{j_1} d_{j_1}(X_{1j_1},x_{j_1}))}\leq \max_{j\in\bar J_\c}(k_j^{-\omega_j})$. Further note that \\$E[k_{j_1}^{\omega_{j_1} d_{j_1}(X_{1j_1},x_{j_1}))}K_\bo(\bX_1,\bx)d_\c(\bbX_1,\bbx)]\leq \max d_\c E[\prod_{j=1}^{p_\text{fun}}K_j(\omega_jd_j(X_{1j},x_j))]$  with $\max d_\c$ being the maximum value that $d_\c$ can take, and that
\beq
E[K_\bo(\bX,\bx)]&=&EE[[K_\bo(\bX,\bx)|\bX_\f]]\\
&=&E\Big[\prod_{j=1}^{p_\f}K_j(\omega_j d_j(X_j,x_j))\cdot E\big[\prod_{j=p_\f+1}^{p_\f+p_\c}(k_j^{-\omega_j})^{d_j(X_j,x_j)}|\bX_\f\big]\Big]\\
&\geq&E\Big[\prod_{j=1}^{p_\text{fun}}K_j(\omega_jd_j(X_{j},x_j))\cdot P(d_\c(\bX,\bx)=0|\bX_\f)\Big]\\
&\geq&E\Big[\prod_{j=1}^{p_\text{fun}}K_j(\omega_jd_j(X_{j},x_j))\Big]\cdot C_\c
\eeq
with Assumption (A0). The same holds with $m$ replaced by $P_g$.

\boxi

%%%%%%%%%%%%%%%%%%%%%%%%%%%%%%%%%%%%%%%

\begin{lemma}\label{EY}
Under Assumptions (A3), (A4), (A5) and (A6) it holds
\[\sup_{\bx\in\cS}\Big|\frac{\frac1{n}\sum_{s=1}^n Z_sK_\bo(\bX_s,\bx)}{E[K_\bo(\bX_1,\bx)]}-\frac{E\big[\frac1{n}\sum_{s=1}^n Z_sK_\bo(\bX_s,\bx)\big]}{E[K_\bo(\bX_1,\bx)]}\Big|=O_P\left(\sqrt{\frac{\log(n)}{n\tilde\phi}}\right)\]
for $Z=Y$ or $Z=I\{Y=g\}$ and all $g=1,\ldots,G$ respectively with $\tilde \phi=\frac{\phi_2^2}{\phi_1}$.
\end{lemma}
The proof of Lemma \ref{EY} follows along the same lines as the proof of Lemma \ref{randden} and is postponed to the online supplement.

%%%%%%%%%%%%%%%%%%%%%%%%%%%%%%%%%%%%%%%%%%%%

{\bf Proof of Theorem \ref{rate}:}
The proof is provided for the classification case. In the regression case the same holds with $I\{Y_i=g\}$ replaced by $Y_i$, $\hat P_{g(-i)}$ by $\hat m_{(-i)}$ and $P_g$ by $m$.

The function $Q$ can be decomposed in
\beq
Q(\omega_1,\ldots,\omega_p)&=&\frac 1n\sum_{i=1}^n\sum_{g=1}^G(I\{Y_i=g\}-\hat P_{g(-i)}(\bX_i))^2v(\bX_i)\\
&=&\frac 1n\sum_{i=1}^n\sum_{g=1}^G (P_g(\bbX_i)-\hat P_{g(-i)}(\bX_i)+I\{Y_i=g\}-P_g(\bbX_i))^2v(\bX_i)\\
&=&\frac 1n\sum_{i=1}^n\sum_{g=1}^G v(\bX_i)(P_g(\bbX_i)-\hat P_{g(-i)}(\bX_i))^2\\
&&+\frac 2n\sum_{i=1}^n\sum_{g=1}^G (P_g(\bbX_i)-\hat P_{g(-i)}(\bX_i)(I\{Y_i=g\}-P_g(\bbX_i))v(\bX_i)\\ 
&&+\frac 1n\sum_{i=1}^n\sum_{g=1}^G (I\{Y_i=g\}-P_g(\bbX_i))^2v(\bX_i).
\eeq
The last term of the decomposition does not depend on $\bo$. In Lemma \ref{epsm} it is shown that 
\[\frac 2n\sum_{i=1}^n\sum_{g=1}^G v(\bX_i)(P_g(\bbX_i)-\hat P_{g(-i)}(\bX_i)(I\{Y_i=g\}-P_g(\bbX_i))=o_P(\zeta(\bbo)+(n\tilde\phi)^{-1})\]
and
\[\frac 2n\sum_{i=1}^n\sum_{g=1}^G v(\bX_i)(P_g(\bbX_i)-\hat P_{g(-i)}(\bX_i)(I\{Y_i=g\}-P_g(\bbX_i))=o_P(1)\ \ \text{uniformly in $\bo$.}\]
Thus, it remains to examine the first term of the decomposition of $Q$. This is
\beq
&&\frac 1n\sum_{i=1}^n\sum_{g=1}^G (\hat P_{g(-i)}(\bX_i)-P_g(\bbX_i))^2\\
%&=&\frac 1n\sum_{i=1}^n \Bigg(\frac{\frac1{n-1}\sum_{s\neq i} I\{Y_s=g\} K(\omega_1d_1(X_{s1},X_{i1}) + \ldots + \omega_pd_p(X_{sp},X_{ip}))}{\frac1{n-1}\sum_{s\neq i} K(\omega_1d_1(X_{s1},X_{i1}) + \ldots + \omega_pd_p(X_{sp},X_{ip}))}-P_g(\bbX_i)\Bigg)^2\\
%&=&\frac 1n\sum_{i=1}^n\sum_{g=1}^G \Bigg(\frac{\frac1{n-1}\sum_{s\neq i}(P_g(\bbX_s)-P_g(\bbX_i))K_{si}}{\frac1{n-1}\sum_{s\neq i}K_{si}} +\frac{\frac1{n-1}\sum_{s\neq i}(I\{Y_s=g\}-P_g(\bbX_s))K_{si}}{\frac1{n-1}\sum_{s\neq i}K_{si}}\Bigg)^2\\
&=&\frac 1n\sum_{i=1}^n v(\bX_i) \Bigg(\frac{\frac1{n-1}\sum_{s\neq i} (P_g(\bbX_s)-P_g(\bbX_i)) K_{si}}{\frac1{n-1}\sum_{s\neq i} K_{si}}+\frac{\frac1{n-1}\sum_{s\neq i} (I\{Y_s=g\}-P_g(\bbX_s))K_{si}}{\frac1{n-1}\sum_{s\neq i} K_{si}}\Bigg)^2\\
&=&\frac 1n\sum_{i=1}^n\sum_{g=1}^G v(\bX_i)\frac{\left(\frac1{n-1}\sum_{s\neq i}(P_g(\bbX_s)-P_g(\bbX_i))K_{si}\right)^2}{\left(\frac1{n-1}\sum_{s\neq i}K_{si}\right)^2}\\
&&+\frac 1n\sum_{i=1}^n\sum_{g=1}^Gv(\bX_i)\frac{\left(\frac1{n-1}\sum_{s\neq i}(I\{Y_s=g\}-P_g(\bbX_s))K_{si}\right)^2}{\left(\frac1{n-1}\sum_{s\neq i}K_{si}\right)^2}\\
&&+\frac 2n\sum_{i=1}^n\sum_{g=1}^Gv(\bX_i)\frac{\left(\frac1{n-1}\sum_{s\neq i}(P_g(\bbX_s)-P_g(\bbX_i))K_{si}\right)\left(\frac1{n-1}\sum_{s\neq i}(I\{Y_s=g\}-P_g(\bbX_s))K_{si}\right)}{\left(\frac1{n-1}\sum_{s\neq i}K_{si}\right)^2}
\eeq
with $K_{si}:=K_\bo(\bX_s,\bX_i)$. 
In Lemma \ref{eps2}  it is shown that
\beq
&&\frac 1n\sum_{i=1}^n\sum_{g=1}^Gv(\bX_i)\frac{\left(\frac1{n-1}\sum_{s\neq i}(I\{Y_s=g\}-P_g(\bbX_s))K_{si}\right)^2}{\left(\frac1{n-1}\sum_{s\neq i}K_{si}\right)^2}\\
&=&\frac {\sum_{g=1}^GE[(I\{Y_1=g\}-P_g(\bbX_1))^2]}{n-1} E\Bigg[\frac{E\big[\big( K_{\bo}(\bX_1,\bX_2)\big)^2|\bX_2\big]}{\Big(E\big[ K_{\bo}(\bX_1,\bX_2)|\bX_2\big]\Big)^2}v(\bX_2)\Bigg]+o_P((n\tilde\phi)^{-1})=o_P(1)\eeq
and in Lemma  \ref{meps} that
\beq&&\frac 2n\sum_{i=1}^n\sum_{g=1}^Gv(\bX_i)\frac{\left(\frac1{n-1}\sum_{s\neq i}(P_g(\bbX_s)-P_g(\bbX_i))K_{si}\right)\left(\frac1{n-1}\sum_{s\neq i}(I\{Y_s=g\}-P_g(\bbX_s))K_{si}\right)}{\left(\frac1{n-1}\sum_{s\neq i}K_{si}\right)^2}\\&=&o_P(\zeta(\bbo))\eeq
and
\beq&&\frac 2n\sum_{i=1}^n\sum_{g=1}^Gv(\bX_i)\frac{\left(\frac1{n-1}\sum_{s\neq i}(P_g(\bbX_s)-P_g(\bbX_i))K_{si}\right)\left(\frac1{n-1}\sum_{s\neq i}(I\{Y_s=g\}-P_g(\bbX_s))K_{si}\right)}{\left(\frac1{n-1}\sum_{s\neq i}K_{si}\right)^2}\\&=&o_P(1)\ \ \text{uniformly in $\bo$.}\eeq
For the remaining term it holds
\beq
&&\frac 1n\sum_{i=1}^n\sum_{g=1}^Gv(\bX_i) \frac{\left(\frac1{n-1}\sum_{s\neq i}(P_g(\bbX_s)-P_g(\bbX_i))K_{si}\right)^2}{\left(\frac1{n-1}\sum_{s\neq i}K_{si}\right)^2}\\
&=&\frac 1n\sum_{i=1}^n\sum_{g=1}^Gv(\bX_i) \frac{\frac1{(n-1)^2}\sum_{s\neq i}(P_g(\bbX_s)-P_g(\bbX_i))^2K_{si}^2}{\left(\frac1{n-1}\sum_{s\neq i}K_{si}\right)^2}\\
&&+\frac 1n\sum_{i=1}^n\sum_{g=1}^Gv(\bX_i) \frac{\frac1{(n-1)^2}\sum_{s_1\neq i}\sum_{\substack{s_2\neq s_1,\\ s_2\neq i}}(P_g(\bbX_{s_1})-P_g(\bbX_i))(P_g(\bbX_{s_2})-P_g(\bbX_i))K_{s_1i}K_{s_2i}}{\left(\frac1{n-1}\sum_{s\neq i}K_{si}\right)^2}.
\eeq
In Lemma \ref{m2} it is shown that
\[\frac 1n\sum_{i=1}^n\sum_{g=1}^Gv(\bX_i) \frac{\frac1{(n-1)^2}\sum_{s\neq i}(P_g(\bbX_s)-P_g(\bbX_i))^2K_{si}^2}{\left(\frac1{n-1}\sum_{s\neq i}K_{si}\right)^2}=o_P(\zeta(\bbo))\]
and
\[\frac 1n\sum_{i=1}^n\sum_{g=1}^Gv(\bX_i) \frac{\frac1{(n-1)^2}\sum_{s\neq i}(P_g(\bbX_s)-P_g(\bbX_i))^2K_{si}^2}{\left(\frac1{n-1}\sum_{s\neq i}K_{si}\right)^2}=o_P(1)\ \ \text{uniformly in $\bo$.}\]
Finally, Lemma \ref{mm} yields for the remaining term
\beq
&&\frac 1n\sum_{i=1}^n\sum_{g=1}^Gv(\bX_i) \frac{\frac1{(n-1)^2}\sum_{s_1\neq i}\sum_{\substack{s_2\neq s_1,\\ s_2\neq i}}(P_g(\bbX_{s_1})-P_g(\bbX_i))(P_g(\bbX_{s_2})-P_g(\bbX_i))K_{s_1i}K_{s_2i}}{\left(\frac1{n-1}\sum_{s\neq i}K_{si}\right)^2}\\
&=&\frac{n-2}{n-1}\sum_{g=1}^GE\Big[\Big(\frac{E\big[P_g(\bbX_1) \bar K_{\bbo}(\bbX_1,\bbX_2)|\bbX_2\big]}{E\big[ \bar K_{\bbo}(\bbX_1,\bbX_2)|\bbX_2\big]}-P_g(\bbX_2)\Big)^2v(\bX_2)\Big]+o_P\big(\zeta(\bbo)\big),
\eeq
and
\beq
&&\frac 1n\sum_{i=1}^n\sum_{g=1}^Gv(\bX_i) \frac{\frac1{(n-1)^2}\sum_{s_1\neq i}\sum_{\substack{s_2\neq s_1,\\ s_2\neq i}}(P_g(\bbX_{s_1})-P_g(\bbX_i))(P_g(\bbX_{s_2})-P_g(\bbX_i))K_{s_1i}K_{s_2i}}{\left(\frac1{n-1}\sum_{s\neq i}K_{si}\right)^2}\\
&=&\frac{n-2}{n-1}\sum_{g=1}^GE\Big[\Big(\frac{E\big[P_g(\bbX_1) \bar K_{\bbo}(\bbX_1,\bbX_2)|\bbX_2\big]}{E\big[ \bar K_{\bbo}(\bbX_1,\bbX_2)|\bbX_2\big]}-P_g(\bbX_2)\Big)^2v(\bX_2)\Big]+o_P(1)\ \ \text{uniformly in $\bo$.}
\eeq
\boxi

%%%%%%%%%%%%%%%%%%%%%%%%%%%%%%%%%%%%%%%%%%%%%

The proofs of the following Lemmas \ref{epsm} - \ref{mm} can be found in the online supplement.  

\begin{lemma}\label{epsm}
Under Assumptions (A0), (A2)--(A7) it holds
 \[\frac 2n \sum_{i=1}^n v(\bX_i)(m(\bbX_i)-\hat m_{-i}(\bX_i))\e_i=O_P(n^{-\frac 12}\sqrt{\zeta(\bbo)})+o_P((n\tilde\phi)^{-1})=o_P(\zeta(\bbo)+(n\tilde\phi)^{-1})\]
and that
\beq
\frac 2n\sum_{i=1}^n\sum_{g=1}^G v(\bX_i)(P_g(\bbX_i)-\hat P_{g(-i)}(\bX_i)(I\{Y_i=g\}-P_g(\bbX_i))&=&O_P(n^{-\frac 12}\sqrt{\zeta(\bbo)})+o_P((n\tilde\phi)^{-1})\\&=&o_P(\zeta(\bbo)+(n\tilde\phi)^{-1})\eeq
with $\zeta(\bbo)=((\min_{j\in\bar J_\f}\omega_j)^{-\beta}+\max_{j\in\bar J_\c}(k_j^{-\omega_j}))^2$ and $\tilde\phi=\frac{\phi_2^2}{\phi_1}$. Note that also $O_P(n^{-\frac 12}\sqrt{\zeta(\bbo)})+o_P((n\tilde\phi)^{-1})=o_P(1)$.
\end{lemma}
%

%%%%%%%%%%%%%%%%%%%%%%%%%%%%%%%%%%%%%%%

\begin{lemma}\label{eps2}
Under Assumptions (A1), (A3), (A4), (A5), (A6) and (A7) it holds
\beq
&&\frac 1n\sum_{i=1}^n\frac{\left(\frac1{n-1}\sum_{s\neq i} \e_s K_\bo(\bX_s,\bX_i)\right)^2}{\left(\frac1{n-1}\sum_{s\neq i} K_\bo(\bX_s,\bX_i)\right)^2}v(\bX_i)\\
&=&\frac {E[\e_1^2]}{n-1} E\Bigg[\frac{E\big[\big( K_{\bo}(\bX_1,\bX_2)\big)^2|\bX_2\big]}{\Big(E\big[ K_{\bo}(\bX_1,\bX_2)|\bX_2\big]\Big)^2}v(\bX_2)\Bigg]+o_P((n\tilde\phi)^{-1})\\
&=&o_P(1),
\eeq
and analogously for the classification case
\beq
&&\frac 1n\sum_{i=1}^n\sum_{g=1}^G\frac{\left(\frac1{n-1}\sum_{s\neq i}(I\{Y_s=g\}-P_g(\bbX_s))K_\bo(\bX_s,\bX_i)\right)^2}{\left(\frac1{n-1}\sum_{s\neq i}K_\bo(\bX_s,\bX_i)\right)^2}v(\bX_i)\\
&=&\frac {\sum_{g=1}^GE[(I\{Y_1=g\}-P_g(\bbX_1))^2]}{n-1} E\Bigg[\frac{E\big[\big( K_{\bo}(\bX_1,\bX_2)\big)^2|\bX_2\big]}{\Big(E\big[ K_{\bo}(\bX_1,\bX_2)|\bX_2\big]\Big)^2}v(\bX_2)\Bigg]+o_P((n\tilde\phi)^{-1})\\
&=&o_P(1),
\eeq
where the $o_P$ terms are uniformly in $\bo$ and $\tilde\phi=\frac{\phi_2^2}{\phi_1}$. \\
Further, the term $E\Big[v(\bX_2)E\big[\big( K_{\bo}(\bX_1,\bX_2)\big)^2|\bX_2\big]/\Big(E\big[ K_{\bo}(\bX_1,\bX_2)|\bX_2\big]\Big)^2\Big]$ is minimized by $\omega_j=0$ for all $j\in\tilde J$.
\end{lemma}
The detailed proof of Lemma \ref{eps2} is postponed to the online supplement. However, the considerations regarding the last assertion of the lemma are presented here.
Consider the term $E\Big[v(\bX_2)E\big[\big( K_{\bo}(\bX_1,\bX_2)\big)^2|\bX_2\big]/\Big(E\big[ K_{\bo}(\bX_1,\bX_2)|\bX_2\big]\Big)^2\Big]$. It holds
\beq
&&E\Bigg[\frac{E\big[\big( K_{\bo}(\bX_1,\bX_2)\big)^2|\bX_2\big]}{\Big(E\big[ K_{\bo}(\bX_1,\bX_2)|\bX_2\big]\Big)^2}v(\bX_2)\Bigg]\\
&=&E\Bigg[\frac{E\big[\big( \bar K_{\bar\omega}(\bbX_1,\bbX_2)\big)^2|\bbX_2\big]}{\Big(E\big[ \bar K_{\bar\omega}(\bbX_1,\bbX_2)|\bbX_2\big]\Big)^2}\frac{E\big[\big( \tilde K_{\tilde\omega}(\tilde\bX_1,\tilde\bX_2)\big)^2|\tilde\bX_2\big]}{\Big(E\big[ \tilde K_{\tilde\omega}(\tilde\bX_1,\tilde\bX_2)|\tilde\bX_2\big]\Big)^2}v(\bX_2)\Bigg],%wegen Unabhängigkeit
\eeq
since $\bar\bX$ and $\tilde\bX$ are independent by Assumption (A1). Now note that
\[0\leq\Var\big(\tilde K_{\tilde\omega}(\tilde\bX_1,\tilde\bX_2)|\tilde \bX_2\big)=E\big[\big(\tilde K_{\tilde\omega}(\tilde\bX_1,\tilde\bX_2)\big)^2|\tilde\bX_2\big]-\Big(E\big[\tilde K_{\tilde\omega}(\tilde\bX_1,\tilde\bX_2)|\tilde\bX_2\big]\Big)^2.\]
If $\omega_j=0$ for all $j\in\tilde J$, it holds that $\tilde K_{\tilde\omega}(\cdot,\cdot)$ is constant and thus  $\Var\big(\tilde K_{\tilde\omega}(\tilde\bX_1,\tilde\bX_2)|\tilde \bX_2\big)=0$ which implies $E\big[\big(\tilde K_{\tilde\omega}(\tilde\bX_1,\tilde\bX_2)\big)^2|\tilde\bX_2\big]/\big(E\big[\tilde K_{\tilde\omega}(\tilde\bX_1,\tilde\bX_2)|\tilde\bX_2\big]\big)^2=1$, so the minimum is achieved. Otherwise, if at least for one $j\in\tilde J$ the weight $\omega_j$ differs from zero, $\Var\big(\tilde K_{\tilde\omega}(\tilde\bX_s,\cdot)\big)>0$ by Assumption (A3) that $K_j$ is not constant and $k_j>1$ for all $j$. Thus, the minimum of $E\big[\big(\tilde K_{\tilde\omega}(\tilde\bX_1,\tilde\bX_2)\big)^2|\tilde\bX_2\big]/\big(E\big[\tilde K_{\tilde\omega}(\tilde\bX_1,\tilde\bX_2)|\tilde\bX_2\big]\big)^2$ is achieved if and only if $\omega_j=0\ \forall j\in\tilde J$.

%%%%%%%%%%%%%%%%%%%%%%%%%%%%%%%%%%%%

\begin{lemma}\label{meps}
Under Assumptions (A0), (A2)--(A7) it holds
\beq
\frac 2n\sum_{i=1}^n v(\bX_i)\frac{\left(\frac1{n-1}\sum_{s\neq i} (m(\bbX_s)-m(\bbX_i)) K_{si}\right)\left(\frac1{n-1}\sum_{s\neq i} \e_s K_{si}\right)}{\left(\frac1{n-1}\sum_{s\neq i} K_{si}\right)^2}&=&O_P(n^{-\frac 12}\sqrt{\zeta(\bbo)})\\&=&o_P(\zeta(\bbo))\eeq
as well as
\beq &&\frac 2n\sum_{i=1}^n\sum_{g=1}^Gv(\bX_i)\frac{\left(\frac1{n-1}\sum_{s\neq i}(P_g(\bbX_s)-P_g(\bbX_i))K_{si}\right)\left(\frac1{n-1}\sum_{s\neq i}(I\{Y_s=g\}-P_g(\bbX_s))\right)}{\left(\frac1{n-1}\sum_{s\neq i}K_{si}\right)^2}\\&=&O_P(n^{-\frac 12}\sqrt{\zeta(\bbo)})=o_P(\zeta(\bbo))\eeq
with $\zeta(\bbo)=((\min_{j\in\bar J_\f}\omega_j)^{-\beta}+\max_{j\in\bar J_\c}(k_j^{-\omega_j}))^2$. Note that also $O_P(n^{-\frac 12}\sqrt{\zeta(\bbo)})=o_P(1)$.
\end{lemma}

%%%%%%%%%%%%%%%%%%%%%%%%%%%%%%%%%%%%%%%

\begin{lemma}\label{m2}
Under Assumptions (A0), (A2), (A3), (A4), (A5) and (A7) it holds
 \[\frac 1n\sum_{i=1}^n v(\bX_i)\frac{\frac1{(n-1)^2}\sum_{s\neq i}(m(\bbX_s)-m(\bbX_i))^2K_{si}^2}{\left(\frac1{n-1}\sum_{s\neq i}K_{si}\right)^2}=O_P(n^{-1}\tilde\phi^{-1}\zeta(\bbo))=o_P(\zeta(\bbo))\]
and accordingly
\[\frac 1n\sum_{i=1}^n\sum_{g=1}^Gv(\bX_i) \frac{\frac1{(n-1)^2}\sum_{s\neq i}(P_g(\bbX_s)-P_g(\bbX_i))^2K_{si}^2}{\left(\frac1{n-1}\sum_{s\neq i}K_{si}\right)^2}=O_P(n^{-1}\tilde\phi^{-1}\zeta(\bbo))=o_P(\zeta(\bbo))\]
with $\zeta(\bbo)=((\min_{j\in\bar J_\f}\omega_j)^{-\beta}+\max_{j\in\bar J_\c}(k_j^{-\omega_j}))^2$. Note that also $O_P(n^{-1}\tilde\phi^{-1}\zeta(\bbo))=o_P(1)$.
\end{lemma}

%%%%%%%%%%%%%%%%%%%%%%%%%%%%%%%%%%%%%%%%%%
\begin{lemma}\label{mm}
Under Assumptions (A0), (A1), (A2), (A3), (A4), (A5) and (A7) it holds
\beq
&&\frac 1n\sum_{i=1}^nv(\bX_i) \frac{\frac1{(n-1)^2}\sum_{s_1\neq i}\sum_{\substack{s_2\neq s_1,\\ s_2\neq i}}(m(\bbX_{s_1})-m(\bbX_i))(m(\bbX_{s_2})-m(\bbX_i))K_{s_1i}K_{s_2i}}{\left(\frac1{n-1}\sum_{s\neq i}K_{si}\right)^2}\\
&=&\frac{n-2}{n-1}E\Big[\Big(\frac{E\big[m(\bbX_1) \bar K_{\bbo}(\bbX_1,\bbX_2)|\bbX_2\big]}{E\big[ \bar K_{\bbo}(\bbX_1,\bbX_2)|\bbX_2\big]}-m(\bbX_2)\Big)^2v(\bX_2)\Big]+O_P\Big(\zeta(\bbo)\phi_2^{-1}\sqrt{\frac{\log(n)}{n\tilde\phi}}\Big)
\eeq
as well as
\beq
&&\frac 1n\sum_{i=1}^n\sum_{g=1}^Gv(\bX_i) \frac{\frac1{(n-1)^2}\sum_{s_1\neq i}\sum_{\substack{s_2\neq s_1,\\ s_2\neq i}}(P_g(\bbX_{s_1})-P_g(\bbX_i))(P_g(\bbX_{s_2})-P_g(\bbX_i))K_{s_1i}K_{s_2i}}{\left(\frac1{n-1}\sum_{s\neq i}K_{si}\right)^2}\\
&=&\frac{n-2}{n-1}\sum_{g=1}^GE\Big[\Big(\frac{E\big[P_g(\bbX_1) \bar K_{\bbo}(\bbX_1,\bbX_2)|\bbX_2\big]}{E\big[ \bar K_{\bbo}(\bbX_1,\bbX_2)|\bbX_2\big]}-P_g(\bbX_2)\Big)^2v(\bX_2)\Big]+O_P\Big(\zeta(\bbo)\phi_2^{-1}\sqrt{\frac{\log(n)}{n\tilde\phi}}\Big)
\eeq
with $\zeta(\bbo)=((\min_{j\in\bar J_\f}\omega_j)^{-\beta}+\max_{j\in\bar J_\c}(k_j^{-\omega_j}))^2$. Note that $O_P\Big(\zeta(\bbo)\phi_2^{-1}\sqrt{\frac{\log(n)}{n\tilde\phi}}\Big)=o_P(\zeta(\bbo))$ and also $O_P\Big(\zeta(\bbo)\phi_2^{-1}\sqrt{\frac{\log(n)}{n\tilde\phi}}\Big)=o_P(1)$.\\
Further, $E\Big[\Big(\frac{E\big[m(\bbX_1) \bar K_{\bbo}(\bbX_1,\bbX_2)|\bbX_2\big]}{E\big[ \bar K_{\bbo}(\bbX_1,\bbX_2)|\bbX_2\big]}-m(\bbX_2)\Big)^2v(\bX_2)\Big]=O(\zeta(\bbo))$ and the same with $m$ replaced by $P_g$.
\end{lemma}

%%%%%%%%%%%%%%%%%%%%%%%%%%%%%%%%%%%%%%%%%%%%%%%%

\section{Concluding remarks}\label{outlook}
The paper revisits the nonparametric method for classification and regression estimation, where the covariates may be functional and categorical, as in \citet{SelkGertheiss2022}.
A thorough theoretical analysis of the asymptotic properties of their method is given. Both the consistency of the regression and classification estimators and the asymptotic behaviour of the weights are derived.

The results of \citet{HallLiRacine2007} are extended to functional combined with categorical covariates. It would also be possible to include continuous variables as well, with very similar arguments. At the same time the results of \citet{FerratyEtal2010} are extended from a single to multiple functional covariates. A further extension is given by the inclusion of the classification case.

%Another topic of interest could be a hypotheses test for significance of the covariates. 
%\citet{RacineHartLi2006} offer such a test for the relevance of categorical covariates in a model with mixed continuous and categorical covariates. However their method cannot be directly adapted to the m

%\bibliographystyle{asa}
\bibliographystyle{apalike}
\bibliography{BibLeonie}

\end{document}